# An Intuitively Complete Analysis of Gödel's Incompleteness

JASON W. STEINMETZ

A detailed and rigorous analysis of Gödel's proof of his first incompleteness theorem is presented. The purpose of this analysis is two-fold. The first is to reveal what Gödel actually proved to provide a clear and solid foundation upon which to base future research. The second is to construct a coherent explication of Gödel's proof that is not only approachable by the non-specialist, but also brings to light the core principles underlying Gödel's proof.

## 1. INTRODUCTION

With Gödel's famous incompleteness proof he attempted to prove some kind of limit on what can be proven within any sufficiently robust and purely formal mathematical system. Subsequently, in the many decades since Gödel published his proof myriad interpretations and theories as well as philosophical extensions of Gödel's proof have been devised. However, we assert that much of this has simply served to confuse and obscure what Gödel actually proved because what a proof actually proves and what a proof is interpreted to prove can be very different and even potentially incommensurate things.

Therefore, the intention of this proof will be to reveal the foundation of Gödel's incompleteness proof within the context of a modern understanding of the fundamental issues that underlie the concept of incompleteness in order to provide a coherent basis for understanding and extending Gödel's ideas. Furthermore, it is hoped that by clearly illustrating the core principles of Gödel's proof that his proof will also become more accessible to those outside the domain of specialists in mathematical logic. However, an introduction to mathematical logic will not be provided and thus a minimal understanding of mathematical logic will be assumed along with at least a rudimentary familiarity with proof theory, model theory, recursive function theory and metamathematics.

Section 2: Background contains an informal introduction to the ideas that form the foundation of Gödel's proof as well as the motivation for the proof. This section is intended to provide a simple intuitive introduction to the concepts and theories along with the formal structures that are the foundation of Gödel's proof. Readers with more than a rudimentary familiarity with Gödel's proof may begin reading at section 2.5: The Fundamental Concepts of Proof and Truth.

In section 3: Notation and Terminology the standard notation from mathematical logic will be formally defined along with some additional notation that is specific to Gödel's proof. In addition, some elementary terminology concerning sets and functions will also be defined.

Section 4: Gödel's Proof contains the formal proof of Gödel's first incompleteness theorem along with a brief description of the proof. Everything that is required to understand the formalization of the proof of Gödel's theorem in this section will have been defined and sufficiently explained in the previous two sections with the exception of Gödel's concept of ω-inconsistency.

Section 5: The Concepts of ω-inconsistency and Mathematical Induction contains the formal definition of Gödel's concept of ω-inconsistency. The customary interpretation of the concept will then also be investigated, which requires an examination of the relationship between the concept of ω-inconsistency and the principle of mathematical induction.

In section 6: The Intuitionistic Acceptability of Gödel's Proof the idea of what may be considered to be a constructive proof will be investigated in an attempt to evaluate Gödel's statement that the proof of his first incompleteness theorem is constructive and thus also intuitionistically acceptable.

Section 7: A Meticulous Intuitive Analysis of Gödel's Proof is a detailed semi-formal examination of the proof of Gödel's first incompleteness theorem.

Section 8: Rosser's Proof contains Kleene's formal proof of Rosser's incompleteness proof from [4] along with a brief description of the proof. Rosser only required the assumption that the system in his proof is simply consistent and thus the customary interpretation of Rosser's proof is that Rosser was able dispense with Gödel's assumption that the formal system is also ω-consistent.

In section 9: <u>The Difference between Gödel's Proof and Rosser's Proof</u> the differences between Gödel's and Rosser's incompleteness proofs are investigated, which includes a new, simple proof of Gödel's theorem that only requires the assumption that the system is simply consistent.

Section 10: <u>Tarski, Mostowski and Robinson's Proof</u> contains the incompleteness proof that is found in *Undecidable Theories* by Tarski in collaboration with Mostowski and Robinson [5].

Finally, in section 11: <u>Summary</u> the concepts of incompleteness and ω-inconsistency are examined in more detail, which includes formal interpretations of ω-inconsistency from Tarski [6.a] and Quine [7.a] along with a simple metamathematical proof that reveals the basic underlying structure of all the incompleteness proofs discussed as well as the role of ω-inconsistency in Gödel's proof.

## 2. BACKGROUND

The background or context within which Gödel published his proof is essential to understanding what Gödel intended to prove and thus also what he actually did prove. Therefore, a brief intuitive description of the underlying question that Gödel attempted to answer with his proof will be provided along with an entirely informal introduction to the fundamental concepts and structures that serve as the foundation of the formal system that Gödel employed within his proof.

### 2.1 Formal Axiomatic Systems

A formal system is, in essence, a system that has been explicitly and completely defined. At its most basic level a formal system consists of a language along with a list of rules that define what constitutes a syntactically well-formed or semantically meaningful formula. The language is comprised of a finite collection of symbols that represent the most primitive elements of the language and are used to construct the formulas of the system. Subsequently, the description of the formal system is slightly different depending on whether the formal system is described from a proof-theoretic or a model-theoretic perspective.

In the most simplistic terms, a proof-theoretic system consists of a list of rules of inference that are used to construct the deductive proofs of the formulas in the system whereas a model-theoretic system consists of a model, which is a collection of objects that the formulas of the system make statements about. Hence, a proof-theoretic system may be characterized as only being concerned with syntax or the actual formulas themselves whereas a model-theoretic system may be characterized as being concerned with semantics or the actual things that the formulas refer to.

An axiomatic system is a system that takes one or more formulas to be the axioms of the system, which may potentially be an infinite number of formulas if an axiom schema is employed. The axioms of the system are a collection of formulas that are asserted to be universally true and from which all the other true formulas or theorems of the system are inferred. In a proof-theoretic system the theorems of the system are deductively proven from the axioms of the system or from previously proven theorems. In a model-theoretic system the axioms of the system define the valid relationships that exist between the objects that constitute the model of the system and thus the theorems of the system are proven based on what is true of the objects within the model.

A proof-theoretic system and a model-theoretic system are not actually two different systems, but are rather two different aspects of any formal axiomatic system. The most significant difference between the two types of systems is how the theorems of the system are proven. A proof-theoretic system is based solely on deductive truth or deductively deriving theorems from the axioms of the system, regardless of what the formulas of the system actually refer to, whereas a model-theoretic system explicitly considers what the formulas of the system refer to and are true of.

Furthermore, a model-theoretic system is typically also a set-theoretic system because the axioms of Set Theory are considered to define what could be called the universal model. Since Set Theory is considered to express the most primitive concept of number then, in the same way that the concept of number is often considered to be the basis of everything, Set Theory is considered to be the basis of every system. This is why Set Theory is sometimes referred to as the theory of everything. The *de facto* standard formalization of Set Theory is Zermelo-Fraenkel Set Theory with the Axiom of Choice or ZFC Set Theory. Therefore, when the specific axioms of Set Theory are directly referenced in what follows, the axioms of ZFC Set Theory are implied.

The formal system that Gödel constructed in his proof is inherently a proof-theoretic system. Gödel does not describe his system as a proof-theoretic system as opposed to a model-theoretic system because the process of differentiating these two types of systems was still in its infancy at the time Gödel published his proof. Thus, this is not a strict designation although it is unequivocal because the defining characteristic of Gödel's system is the concept of a proof-theoretic deductive proof within a formal axiomatic system.

It is of interest to note that a complete system is often defined in modern texts as a system within which a proof-theoretic deductive proof is equivalent to a model-theoretic proof. Therefore, within this context, if a formal system is complete then the existence of a proof in either system implies the existence of a proof in the other.

### 2.2 Consistency, Completeness and Hilbert's Program

The underlying question that Gödel's proof was intended to answer is related to what has become known as Hilbert's program. The primary objective of Hilbert's program was to construct a simply consistent and complete formal axiomatic system that could be proven to be simply consistent and complete using only finitary methods of proof and within which all the essential theorems of mathematics could be proven. Subsequently, since Gödel's incompleteness theorems are generally interpreted to have proven that any sufficiently robust and simply consistent formal axiomatic system must also be incomplete, Gödel's proof is also generally interpreted to have proven that Hilbert's program cannot succeed.

More intuitively, the objective of Hilbert's program was to construct a formal axiomatic system that is unquestionably valid and also capable of a proving whether any meaningful sentence in the language of the system is true. Since a formal axiomatic system is defined to be simply consistent if a contradiction cannot be inferred from its axioms then the formal requirement that the system be simply consistent requires the system to be incapable of proving any sentence to be both true and false, which is a contradiction. The formal requirement that the system be complete requires the system to be capable of proving every meaningful sentence to be either true or false. Since a sentence that is true may be interpreted to be the refutation of the negation of that sentence then a formal axiomatic system may alternately be defined to be complete if every meaningful sentence that is true can be inferred from its axioms. As a result, the requirement that the formal system be proven to be simply consistent and complete using only finitary methods of proof serves to establish the unquestionable validity and unrestricted applicability of the system. This requirement is informal because the concept of a finitary method of proof has never been precisely defined.

What Hilbert called a finitary method of proof may be characterized as method of proof that is inherently finite and thus does not depend on any questionable assumptions about infinity or assumptions about what is true of a potentially or actually infinite number of objects. Since some of the formulas being considered may apply to an infinite number of objects, it is the actual method that is characterized as finite. Therefore, within the context of Hilbert's program, a finitary method of proof suggests a method of proof that employs only finitary reasoning that can be easily and clearly justified. Hence, Hilbert's emphasis on finitary methods of proof may be characterized as an emphasis on methods of proof that are not controversial and thus are also unquestionably and universally acceptable.

It is of interest to note that in Gödel's second incompleteness theorem Gödel formally defined a simply consistent formal axiomatic system as a system within which at least one sentence cannot be proven to be true regardless of whether the sentence is actually true. This definition of simple consistency is based on *ex falso quodlibet* or the principle of explosion, which asserts that every sentence can be inferred from the axioms of a simply inconsistent system.

The pertinent details of Hilbert's program can be found in van Heijenoort [2], particularly [2.b].

### 2.3 Effective Computability

The concept of effective computability is an important aspect of Gödel's proof. The general characterization of an effectively computable function that Gödel defined in his proof is not only a crucial part of his proof, but it also became the basis of the modern theory of recursive functions. The recursive functions are one of three models of computation that were subsequently proven to be equivalent and, according to Church's Thesis, serve to define what can be effectively computed by any model of computation.

The recursive functions consist of the primitive recursive functions and the general recursive functions. Gödel defined the primitive recursive functions in his proof, although he referred to them simply as recursive because the modern terminology had not yet been established. The difference between the two classes of functions is that, due to how the primitive recursive functions are defined, their computation is bounded whereas the general recursive functions may employ an unbounded search. As a result, the computation of a primitive recursive function is always finite or will always eventually terminate in a finite number of steps. However, the computation of a general recursive function may potentially never terminate. Therefore, since a computation that never terminates is effectively synonymous with a computation that does not compute a value, a general recursive function will not necessarily always compute a value whereas a primitive recursive function will always compute a value.

Gödel defined a list of relations in his proof that precisely define the basic operations that he employed in his proof (Appendix B: Gödel's System). Subsequently, in a parenthetical note after the definition of the last relation, relation 46, Gödel asserts that every relation that he defined can be proven to be [primitive] recursive with the single exception of relation 46, which asserts that a primitive recursive function can be defined to implement or compute every relation except relation 46. The reason relation 46 cannot be asserted to be primitive recursive is because the definition of relation 46 employs an unbounded search whereas every search within the definition of all of the other relations is explicitly bounded.

Therefore, since the relations that Gödel employed to prove his incompleteness theorems are primitive recursive and thus the relations represent inherently finite computations then the role of the concept of effective computability within Gödel's proof may be characterized as fulfilling the requirement of Hilbert's program that the methods of proof be finitary. As a consequence, the concept of effective computability serves to substantiate the validity of Gödel's proof because the methods that Gödel employed in his proof do not involve any subjective processes or steps that cannot be completely and objectively analyzed and validated. Although relation 46 is not primitive recursive and thus that relation does not necessarily represent a computation that is inherently finite, it may still be possible to consider the relation to be finitary due to how its use is restricted within Gödel's proof. This will be specifically addressed in more detail in what follows.

### 2.4 Gödel Numbers and the Arithmetization of Syntax

Gödel constructed what have come to be known as Gödel numbers. A Gödel number is the unique numeric representation of a unique sequence of symbols that represent a unique formula in the language of a formal system. The process employed to construct Gödel numbers is often referred to as the arithmetization of syntax because a unique sequence of symbols, which is the syntax of a formula, is transformed into a unique natural number that can be analyzed and manipulated by a mathematical system based on the concepts of arithmetic. Furthermore, Gödel also proved that every unique finite sequence of formulas in a formal system also corresponds to a unique natural number and thus a Gödel number is also the unique numeric representation of a unique finite sequence of formulas in a formal system.

Each of the primitive symbols in the language of Gödel's formal system is assigned the unique natural number 1, 3, 5, 7, 9, 11, 13 and then every variable is assigned a prime number greater than 13 with an exponent greater than zero indicating its type. Hence, a formula, which is a finite sequence of $n \geq 1$ symbols, is also a sequence of natural numbers. We will denote this sequence by $m_1, \ldots, m_n$, where each $m_i$ is the natural number that was assigned to the $i^{\text{th}}$ symbol in the formula. The Gödel number of a formula is then defined as $p_1^{m_1} \times \ldots \times p_n^{m_n}$, where each $p_i$ denotes the $i^{\text{th}}$ prime number. This number is unique due to the prime decomposition theorem, which proves that every unique natural number has a unique decomposition into its prime components. The Gödel number of a finite sequence of formulas is constructed similarly. A finite sequence of $n \geq 1$ formulas is also a sequence of Gödel numbers. We will denote this sequence by $G_1, \ldots, G_n$, where each $G_i$ is the Gödel number of the $i^{\text{th}}$ formula in the sequence. As a result, the Gödel number of a finite sequence of formulas is $p_1^{G_1} \times \ldots \times p_n^{G_n}$, where again each $p_i$ denotes the $i^{\text{th}}$ prime number.

In the preceding construction the type of a variable is referred to. This is a reference to the simple theory of types, which Gödel employed in his proof. However, the simple theory of types is not directly relevant to our examination of what Gödel proved. Therefore, any reference to the type of a variable in what follows is included for technical accuracy only and may be ignored.

This construction of the Gödel numbers is more intuitively expressive than practical or technically efficient. Gödel employed the prime numbers to construct his Gödel numbers because, in the same way that the prime numbers are intuitively perceived to be the primitive components of a natural number, the symbols of a formula are its primitive components. Therefore, due to how the Gödel numbers are constructed, a Gödel number relates to a formula in the same way that a natural number relates to its prime components. Hence, the construction of Gödel numbers serves to lend credence to the intuitive idea that the formulas in a formal system may be treated as mathematical objects or numeric objects within a mathematical system based on the concepts of arithmetic.

The principal assertion of the arithmetization of syntax is that the Gödel number of every formula and every finite sequence of formulas in a formal system exists and is a unique natural number. As a result, since every deductive proof in a formal system is some finite sequence of formulas then this allowed Gödel to treat the formulas and the deductive proofs in a formal system as numbers in a mathematical system. Furthermore, since the arithmetization of syntax also permits the complete enumeration of every formula and every deductive proof in a formal system in the same way that every natural number may be completely enumerated then the arithmetization of syntax also allowed Gödel to quantify and thus reason over the totality of every formula and every deductive proof in a formal system.

A more modern exposition of the arithmetization of syntax that includes a formal definition of the one-to-one correspondence between the well-formed formulas in a formal system and the natural numbers is provided in Appendix A: Gödel Numbers and the Arithmetization of Syntax.

### 2.5 The Fundamental Concepts of Proof and Truth

Prior to the publication of Gödel's proof the discussion about mathematical proofs in the literature generally consisted of whether a formula was true or not and whether the truth of that formula could be conclusively proven within a formal mathematical system. This is an inherently naïve and uncritical view of the concepts of proof and truth. The naïve view of the concepts of proof and truth explicitly assumes that a formal proof is a representation of truth. However, this is an invalid objectification of the concept of truth. A formal proof is objective, by definition, whereas the concept of truth is not due to its inherently informal nature. As a result, what can be proven within any formal system and what may be considered to be actually true are fundamentally different and even potentially incommensurate ideas.

The primary focus of Gödel's proof was the proof-theoretic concept of proof or, more specifically, what can be deductively proven within any formal axiomatic system that is based on the principles of elementary arithmetic and mathematical logic. The concepts of proof and truth and how they are related is crucial to understanding Gödel's proof because the exclusive focus of Gödel's proof was to determine what can be deductively inferred from the axioms of a formal axiomatic system without any explicit regard for what those axioms may or may not be true of and thus, by extension, without any explicit regard for what may or may not actually be true.

Since a concept of truth may only be circularly defined as what is true, a naïve and uncritical view of Gödel's proof naturally leads to many speculative and unjustified interpretations of what Gödel actually proved. Therefore, in an attempt to avoid any such speculative assumptions we will simply define three propositions that attempt to characterize the truth of proof within any formal system and are, for the most part, well known and not controversial. These propositions are explicitly stated because, regardless of their seemingly obvious validity, they are often violated when speaking about what is true in mathematics.

Within the context of a formal system the concept of truth is synonymous with what can proven within the system.

> **The Truth of Proof, Proposition 1:** within any formal system it is never valid to speak about truth or what is, in reality, true or false, but only about what can and cannot be proven.

Although this proposition implicitly relates the concepts of truth and reality, no such relationship is actually asserted to exist. Some philosophies of mathematics propose that mathematical truth is not dependent on being practical or realistic and this proposition is not intended to refute that perspective. Hence, the idea that is expressed within the proposition by the term reality must be interpreted as being just as subjective and open to interpretation as the idea expressed by the term truth. What is implicitly suggested by contrasting these ideas is more important than whether any explicit and definable relationship actually exists between them.

Within the context of a formal axiomatic system a proof only represents what can be inferred from the axioms of the system.

**The Truth of Proof, Proposition 2:** the truth of a proof within any formal axiomatic system is absolutely relative to the scope and validity of the axioms of the system.

The term scope refers to the universe that the axioms are applicable to and the term validity refers to the extent of that applicability or how correct the axioms are within that universe. We use the term absolutely relative to imply that the scope and validity of any formal axiomatic system is inextricably dependant on and thus also constrained by the scope and validity of its axioms.

A formal axiomatic system is merely the formal analogue of the intuitive idea that is basis of the system. We call this intuitive idea the naïve theory.

**The Truth of Proof, Proposition 3:** the naïve theory of a formal axiomatic system, which the axioms of the system attempt to formally characterize, defines truth within the system.

In essence, this proposition is merely a reformulation of Tarski's characterization of truth within a formal system ([6] pg. 273). Since metamathematics is customarily defined as the language that is employed to speak about formal mathematical languages and formal mathematical systems then metamathematics may be considered to be the language of the naïve theory of mathematics and mathematical logic. Hence, this proposition states that truth is defined within a formal axiomatic system by the metamathematical statements and methods that have been accepted as true, provided the analogous axioms and methods in the formal system are accurate formal characterizations of those metamathematical statements and methods.

As a consequence, although we will be primarily focused on the underlying truth of Gödel's incompleteness proof, it must be kept in mind that Gödel's proof is actually concerned with what can be deductively proven within any sufficiently robust formal axiomatic system based on the principles of elementary arithmetic and mathematical logic.

## 3. NOTATION AND TERMINOLOGY

The notation and terminology that constitutes the basic framework of Gödel's system will now be defined. Only the notation that is required to prove Gödel's first incompleteness theorem along with some additional notation and terminology that is needed to explain Gödel's Proof will be defined. Gödel's full system is reproduced with updated notation and briefly described in Appendix B: Gödel's System.

The notation that we will employ to state a definition is $\equiv_{df}$. Therefore, "$x \equiv_{df} y$" states that "$x$ is defined as $y$" or "the definition of $x$ is $y$." In the following definitions the variables $x$ and $y$ are employed as general, undefined variables.

We will begin by defining the notation for the standard logical operations:

**Negation (Logical NOT):** $\neg x \equiv_{df}$ the negation of $x$.

**Conjunction (Logical AND):** $x \wedge y \equiv_{df}$ the conjunction of $x$ and $y$.

**Disjunction (Logical OR):** $x \vee y \equiv_{df}$ the disjunction of $x$ and $y$.

**Implication (Logical Conditional):** $x \rightarrow y \equiv_{df} \neg x \vee y$

**Equivalence (Logical Biconditional):** $x \leftrightarrow y \equiv_{df} (x \rightarrow y) \wedge (y \rightarrow x)$

We will also define the standard quantifiers from predicate logic:

**Existential Quantification:** $\exists x \equiv_{df}$ there exists some $x$.

**Universal Quantification:** $\forall x \equiv_{df}$ for all $x$.

Although the preceding operations are the basic operations of mathematical logic and thus are found in any introduction to mathematical logic, the notation is explicitly defined because it is not always standard and was different in Gödel's proof.

It should be explicitly noted that often either conjunction or disjunction is taken as primitive, but not both as we have done, and then the other is defined in terms of the first. Existential and universal quantification are also often similarly defined in terms of each other as Gödel did within his formal system.

We will also explicitly define the notation for set membership, which is the primitive operation of Set Theory, although the notation is standard in all modern texts.

**Set Membership:** $x \in Y \equiv_{df} x$ is a member of or contained in the set $Y$.

**Negation of Set Membership:** $x \notin Y \equiv_{df} \neg(x \in Y)$

We will use other standard mathematical notation without definition because it is assumed that the operations and their notation are well known.

We will now briefly define some of the basic concepts of functions within a set-theoretic context, all of which are standard in modern texts.

**Function:** $f: X \rightarrow Y \equiv_{df}$ the function $f$ is a function from the set $X$ to the set $Y$.

The definition simply states that the function $f$ takes arguments from the set $X$ and returns values from the set $Y$. If a function is defined for every argument or every $x \in X$ then it is called a total function, otherwise it is called a partial function. However, the concept of a partial function will not be needed in what follows and thus will not be explicitly defined. Hence, every function is assumed to be a total function.

In the following definition we will employ the customary notation for multiple quantifiers. Hence, $\forall x_0 \in X\ \forall x_1 \in X$ is stated as $\forall x_0, x_1 \in X$, where $x_0 \in X$ and $x_1 \in X$ denote any members of $X$.

**Injective or 1-to-1 Function:** $\forall x_0, x_1 \in X\ (f(x_0) = f(x_1) \rightarrow x_0 = x_1)$

A function is injective or 1-to-1 if the function never returns the same value for distinct arguments or equivalently the value of the function is a unique $y \in Y$ for every distinct argument $x \in X$.

**Surjective or Onto Function:** $\forall y \in Y\ \exists x \in X\ (f(x) = y)$

A function is surjective or onto if every $y \in Y$ is the value of the function for some argument or equivalently no $y \in Y$ exists such that $f(x) \neq y$ for every argument $x \in X$.

**Bijective Function:** $\forall y \in Y\ \exists x_0 \in X\ \forall x_1 \in X\ (f(x_0) = y \wedge \{f(x_1) = y \rightarrow x_0 = x_1\})$

A function is bijective if every $y \in Y$ is the value of the function for some argument $x \in X$ and the value of the function is a unique $y \in Y$ for every distinct argument $x \in X$. A bijective function is a function that is both injective and surjective. A bijection, which is what a bijective function establishes between the sets $X$ and $Y$, is often referred to as a one-to-one correspondence. However, a one-to-one correspondence does not necessarily define a bijective function unless and only unless the set-theoretic definition of a function is assumed along with the Axiom of Choice from Set Theory.

We will also note that if a function $f$ is a bijective function then the inverse of $f$, which is typically denoted by $f^{-1}$, always exists and is also a bijective function. As a result, if $f: X \rightarrow Y$ is a bijective function then its inverse $f^{-1}: Y \rightarrow X$ is also a bijective function. However, the concept of an inverse function will not be needed in what follows and thus it will not be explicitly defined.

We will now define the least search operator from recursive function theory. Although the notation that we employ is due to Kleene ([4] pg. 279) and is the standard notation in most modern texts, it is different from the notation that Gödel employed in his proof.

**Least Search Operator:** $\mu x \equiv_{df}$ the least $x$.

The least search operator is defined in a similar manner to how the quantifiers from predicate logic were defined because the least search operator is effectively identical to the existential quantifier with the sole exception that the existential quantifier is concerned with any possible value of the variable $x$ while the least search operator is specifically concerned with the least or smallest value of the variable $x$ within some ordering of all the possible values of the variable $x$. In general, the least search operator suggests the obvious and rudimentary method of searching for the least value of the variable $x$ by employing a linear search through every possible value of the variable $x$, which is a search that begins by evaluating the first or smallest value of $x$ and then progressively evaluating the next larger value of $x$ until the desired value is found.

Finally, we will define the notation that is not standard or that is specific to the present proof. Throughout the proof we will employ the symbols $\varphi$ and $\Psi$ to represent formulas. If a formula $\varphi$ contains no free variables then $\varphi$ is called a sentence, otherwise for any formula $\varphi$ with one or more free variables:

**Formula with Free Variables:** $\varphi(x_1, \ldots, x_n) \equiv_{\text{df}} \varphi$ has the $n \geq 1$ free variables $x_1, \ldots, x_n$.

This notation may be employed in two different ways, which will be apparent from the context. When $\varphi$ denotes a formula that has been specifically defined then the notation $\varphi(x_1, \ldots, x_n)$ merely indicates the free variables that actually exist in the formula, but when $\varphi$ denotes a formula that has not been specifically defined then the notion $\varphi(x_1, \ldots, x_n)$ asserts that the formula $\varphi$ contains the free variables $x_1, \ldots, x_n$. Also note that each free variable $x_i$ in $\varphi(x_1, \ldots, x_n)$ does not denote a variable with a fixed subscript or a variable that is associated with any specific symbol or symbols, but rather simply denotes the $i^{\text{th}}$ free variable in the formula $\varphi$.

For some formula $\varphi(x_1, \ldots, x_n)$ and some $y_1, \ldots, y_n$, where each $y_i$ denotes a fixed value or is one of the free variables in $\varphi(x_1, \ldots, x_n)$:

**Substitution:** $\varphi(y_1, \ldots, y_n) \equiv_{\text{df}}$ an $m$-tuple of fixed values $c_1, \ldots, c_m$ are substituted for $m$ free variables in $\varphi(x_1, \ldots, x_n)$, where $1 \leq m \leq n$ and each fixed value $c_j$ is some $y_i$ in $\varphi(y_1, \ldots, y_n)$.

Each of the $m \geq 1$ fixed values $c_1, \ldots, c_m$ is some constant or defined term with no free variables. If $m < n$ then each $y_i$ in $\varphi(y_1, \ldots, y_n)$ that is not a fixed value must be the corresponding free variable $x_i$ in $\varphi(x_1, \ldots, x_n)$. The operation of substitution considers a formula to be a sequence of symbols and then, if applicable, substitutes some of those symbols for other symbols and thus constructs a different formula. Substitution is defined in Gödel's formal system by his primitive recursive relation 31, although only for a single free variable (Appendix B: Gödel's System or [1.a] pg. 167, within which Gödel uses the notation $Sb(a\,{}^{x}_{c})$, where $a$ is some $\varphi(x_1, \ldots, x_n)$, $x$ is one of the free variables in $\varphi(x_1, \ldots, x_n)$ and $c$ is some fixed value). Subsequently, Gödel defined the notation for the substitution of multiple free variables in a footnote ([1.a] pg. 167 footnote 37), which is effectively equivalent to the notation that we have just defined.

The present notation was chosen not only because it is more concise and thus more readable, but also because the notation includes more information than the notation that Gödel employed in his proof. Specifically, the notation explicitly indicates the free variables in the formula for which no fixed value is substituted and which are thus still free.

We will now define the notation that represents the process indicated by the arithmetization of syntax (section 2.3: Gödel Numbers and the Arithmetization of Syntax). For any formula $\varphi$:

**Gödel Number of a Formula:** $[\varphi_n] = n$

For every formula $\varphi$ the Gödel number of $\varphi$ exists and is unique, in symbols: $[\varphi] = n$, where $[\varphi]$ indicates the process of constructing the unique natural number $n$ that is the Gödel number of the formula $\varphi$. Hence, we use the notation $\varphi_n$ to denote the formula whose Gödel number is $n$ and then we use $n$ to denote the natural number that is the Gödel number of the formula $\varphi_n$.

Note that for any formula $[\varphi_m(x_1, \ldots, x_n)] = m$, $[\varphi_m(y_1, \ldots, y_n)] \neq m$, where $n \geq 1$ and one or more of the $y_i$ is some fixed value. This is because $m$ is the Gödel number of the formula $\varphi_m$ with the free variables $x_1, \ldots, x_n$ and thus if a fixed value is substituted for one or more of the free variables in $\varphi_m(x_1, \ldots, x_n)$ then a different formula is constructed. As a result, this formula then has a different Gödel number and thus $[\varphi_m(y_1, \ldots, y_n)] \neq m$.

We will now employ the present notation to define Gödel's relations 45 and 46.

**Gödel's Proof Relation:** $\varphi_x$ ProofOf $\varphi_y \equiv_{\text{df}} \varphi_x$ is a finite sequence of formulas that constitutes a deductive proof of $\varphi_y$.

This is Gödel's primitive recursive relation 45 (Appendix B: Gödel's System or [1.a] pg. 171, within which Gödel uses the notation $x\ B\ y$).

**Gödel's Provability Predicate:** $\text{Prov}(\varphi_x) \equiv_{\text{df}} \exists n\ (\varphi_n \text{ ProofOf } \varphi_x)$

This is Gödel's relation 46 (Appendix B: Gödel's System or [1.a] pg. 171, within which Gödel uses the notation $Bew(x)$). As stated in section 2.4: Effective Computability, this relation is not primitive recursive because the existential quantifier in the definition is not explicitly bounded. As a result, if a deductive proof of $\varphi_x$ does not exist within the formal system and thus $\text{Prov}(\varphi_x)$ is false then any linear search for a deductive proof of the formula $\varphi_x$ will never terminate.

Gödel also defined the relations ProofOf and Prov relative to a formal system that was denoted by $\kappa$. Hence, Prov($x$) referred to a proof within Gödel's formal system and $\text{Prov}_\kappa(x)$ referred to a proof within a formal system $\kappa$ that is an extension of Gödel's formal system and thus includes the axioms of Gödel's formal system possibly along with other axioms. Since we will focus only on what Gödel proved within his formal system and not on what might possibly be proven within any extension of Gödel's formal system then this notation will not be explicitly defined. However, it is to be understood that what Gödel proved within his formal system can also be proven within any extension of Gödel's formal system, which Gödel stated in his proof ([1.a] pg. 173, where Gödel asserts $\forall x\ (\text{Prov}(x) \rightarrow \text{Prov}_\kappa(x))$).

## 4. GÖDEL'S PROOF

We will first state Gödel's theorem 5 because the theorem plays a crucial role in his proof and in what follows.

**Gödel's Theorem 5:** for any primitive recursive relation $R(x_1, \ldots, x_n)$, where $n \geq 1$, there exists a primitive recursive numeric formula $\varphi$ with $n$ free variables such that:

1. $R(x_1, \ldots, x_n) \rightarrow \text{Prov}(\varphi(x_1, \ldots, x_n))$
2. $\neg R(x_1, \ldots, x_n) \rightarrow \text{Prov}(\neg\varphi(x_1, \ldots, x_n))$

The outline of the proof of theorem 5 that Gödel provided is what is known in modern texts as a proof by induction on complexity; the complexity of the formula being defined using Gödel's formal definition of primitive recursive functions.

Gödel intuitively characterized theorem 5 as asserting that any primitive recursive relation is definable within his formal system. Due to the definition of a primitive recursive relation, every primitive recursive relation has a primitive recursive function that implements or computes the relation. As a result, the deductive proof of the numeric formula $\varphi(x_1, \ldots, x_n)$ or $\neg\varphi(x_1, \ldots, x_n)$ is simply the explicit calculation of the primitive recursive function that implements or computes the relation $R(x_1, \ldots, x_n)$ within the formal system and this explicit calculation constitutes a deductive proof or refutation of the primitive recursive relation with which it corresponds for the specific numeric values that are substituted for its free variables.

We will now begin Gödel's proof by defining the relation $Q$:

**1:** $Q(x, y) \equiv_{\text{df}} \neg(\varphi_x \text{ ProofOf } \varphi_y(y))$

Due to Gödel's theorem 5 there is a numeric formula $\varphi_q(x, y)$ that implements or computes the relation $Q$. Therefore, due to theorem 5 and the definition of the relation $Q$:

**2:** $\neg(\varphi_x \text{ ProofOf } \varphi_y(y)) \rightarrow \text{Prov}(\varphi_q(x, y))$

**3:** $\varphi_x \text{ ProofOf } \varphi_y(y) \rightarrow \text{Prov}(\neg\varphi_q(x, y))$

The formulas $\varphi_p$ and $\varphi_r$ are now defined:

**4:** $\varphi_p(x) \equiv_{\text{df}} \forall n\ \varphi_q(n, x)$

And then where $p = [\varphi_p(x)]$:

**5:** $\varphi_r(x) \equiv_{\text{df}} \varphi_q(x, p)$

Hence, it is immediate due to 5:

**6:** $\varphi_q(x, p) = \varphi_r(x)$

Substituting $p = [\varphi_p(x)]$ for $x$ in 4 and then due to 6:

**7:** $\varphi_p(p) = \forall n\ \varphi_q(n, p) = \forall n\ \varphi_r(n)$

Finally, substituting $p = [\varphi_p(x)]$ for $y$ in 2 and 3 and then due to 6 and 7:

**8:** $\neg(\varphi_x \text{ ProofOf } \forall n\ \varphi_r(n)) \rightarrow \text{Prov}(\varphi_r(x))$

**9:** $\varphi_x \text{ ProofOf } \forall n\ \varphi_r(n) \rightarrow \text{Prov}(\neg\varphi_r(x))$

First it is assumed that a deductive proof of $\forall n\ \varphi_r(n)$ exists, in symbols: $\exists x\ (\varphi_x \text{ ProofOf } \forall n\ \varphi_r(n))$. Therefore, if a deductive proof of $\forall n\ \varphi_r(n)$ exists then a deductive proof of $\neg\varphi_r(x)$ must also exist for some value of the variable $x$, due to 9, in symbols: $\exists x\ \text{Prov}(\neg\varphi_r(x))$. However, if the formal system is simply consistent then a deductive proof of $\neg\varphi_r(x)$ for some value of the variable $x$ contradicts a deductive proof of $\forall n\ \varphi_r(n)$. Subsequently, it is inferred that it cannot be the case that a deductive proof of $\forall n\ \varphi_r(n)$ exists, in symbols: $\forall x\ \neg(\varphi_x \text{ ProofOf } \forall n\ \varphi_r(n))$.

It is then assumed that a deductive proof of $\neg\forall n\ \varphi_r(n)$ exists, which is a deductive proof of the refutation of $\forall n\ \varphi_r(n)$, and thus a deductive proof of $\forall n\ \varphi_r(n)$ cannot exist if the formal system is assumed to be simply consistent, in symbols: $\forall x\ \neg(\varphi_x \text{ ProofOf } \forall n\ \varphi_r(n))$. This is also the result of the contradiction inferred from 9 as well. Therefore, if a deductive proof of $\forall n\ \varphi_r(n)$ does not exist then a deductive proof of $\varphi_r(x)$ must exist for every value of the variable $x$, due to 8, in symbols: $\forall x\ \text{Prov}(\varphi_r(x))$. However, if a deductive proof of $\varphi_r(x)$ exists for every value of the variable $x$ and a deductive proof of $\neg\forall n\ \varphi_r(n)$ also exists, this is what Gödel defined as an ω-inconsistency. As a result, if the formal system is ω-consistent then a deductive proof of $\varphi_r(x)$ for every value of the variable $x$ contradicts a deductive proof of $\neg\forall n\ \varphi_r(n)$. Subsequently, it is inferred that it cannot be the case that a deductive proof of $\neg\forall n\ \varphi_r(n)$ exists. ■

As a consequence, if the formal system is assumed to be simply consistent then a contradiction results from the assumption that a deductive proof of the sentence $\forall n\ \varphi_r(n)$ exists. And if the formal system is assumed to be ω-consistent then a contradiction results from the assumption that a deductive proof of the sentence $\neg\forall n\ \varphi_r(n)$ exists, which is a deductive proof of the refutation of $\forall n\ \varphi_r(n)$. Therefore, since both assumptions result in a contradiction it is then inferred that a deductive proof of the sentence $\forall n\ \varphi_r(n)$ cannot exist within the formal system and a deductive proof of its refutation also cannot exist within the formal system. Hence, the sentence $\forall n\ \varphi_r(n)$ cannot be proven within the formal system and is thus undecidable within the formal system.

Furthermore, the sentence $\varphi_p(p)$ can easily be seen to be true according to Gödel ([1.a] pg. 151). Since $\varphi_p(p) = \forall n\ \varphi_r(n)$, due to 7, and $\forall n\ \varphi_r(n) = \forall n\ \neg(\varphi_n \text{ ProofOf } \varphi_p(p))$, due to the definition of the formula $\varphi_r$ and the definition of the relation $Q$, then $\forall n\ \varphi_r(n)$ states that a deductive proof of the sentence $\varphi_p(p)$ does not exist and thus $\varphi_p(p)$ states that a deductive proof of itself does not exist, due to 7. However, as was just proven, $\forall n\ \varphi_r(n)$ is an undecidable sentence and thus a deductive proof of the sentence $\forall n\ \varphi_r(n)$ cannot exist, in symbols: $\forall x\ \neg(\varphi_x \text{ ProofOf } \forall n\ \varphi_r(n))$. As a result, the sentence $\varphi_p(p)$ is clearly true, due to 7. Hence, although the sentence $\varphi_p(p)$ can clearly be seen to be true this does not contradict its formal undecidability because, according to Gödel, the truth of the sentence $\varphi_p(p)$ is due to metamathematical considerations.

In conclusion, the Gödel sentence $\varphi_p(p) = \forall n\ \varphi_r(n) = \forall n\ \neg(\varphi_n \text{ ProofOf } \varphi_p(p))$ cannot be proven within the formal system and thus its truth cannot be decided within the formal system. However, its truth is implied by the proof of its undecidability. Therefore, it is asserted that the formal system is inherently incomplete because a legitimate sentence exists in the language of the formal system that can clearly be seen to be true, but which is undecidable within the formal system.

## 5. THE CONCEPTS OF ω-INCONSISTENCY AND MATHEMATICAL INDUCTION

We will begin with Gödel's definition of ω-inconsistency. A theorem is defined as a sentence for which a proof exists within the formal system. The theory of a formal system is then defined as the set that contains every theorem as well as every axiom, in symbols: $Th = \{\varphi \mid \text{Prov}(\varphi)\}$, where $\varphi$ denotes any sentence in the language of the formal system. Intuitively, the theory of a formal system is the set that contains every sentence in the language of the formal system that is true. Subsequently, Gödel defines an ω-inconsistent formal system as, $\forall x\ (\varphi(x) \in Th) \wedge \neg\forall x\ \varphi(x) \in Th$, where $\varphi(x)$ denotes any formula and the values of the variable $x$ range over the natural numbers. This definition states that a deductive proof exists within the formal system for a formula $\varphi(x)$ when any value is substituted for its free variable and a deductive proof also exists within the formal system for the negation of the sentence $\forall x\ \varphi(x)$. As a result, an ω-inconsistency is often defined as a formal system within which a proof exists for every sentence $\varphi(0)$, $\varphi(1)$, $\varphi(2)$, … as well as the sentence $\neg\forall x\ \varphi(x)$, where the values of the variable $x$ range over the natural numbers. The definition of ω-inconsistency that we will use is effectively identical to Gödel's definition, which is, $\forall x\ \text{Prov}(\varphi(x)) \wedge \text{Prov}(\neg\forall x\ \varphi(x))$, where the values of the variable $x$ range over the natural numbers.

Before constructing a precise and fully general definition of ω-inconsistency we will first define a bijective function that enumerates every possible value of the variable $x$. This bijective function is technically defined as: $f$: $\mathbb{N} \rightarrow X$, where $\mathbb{N}$ denotes the set of natural numbers and $X$ denotes the set of every possible value that the variable $x$ can take. We then also define $x_n = f(n)$, where $n \in \mathbb{N}$ and $x_n \in X$. If $X$ is the set of natural numbers or $X = \mathbb{N}$ then $f$ is the identity function or $f(x) = x$ and thus $x_n = f(n) = n$.

Therefore, for any formula $\varphi(x)$ in some formal system and some bijective function $f$: $\mathbb{N} \rightarrow X$ and thus also $x_n = f(n)$:

**ω-inconsistency:** $\forall n \in \mathbb{N}\ \text{Prov}(\varphi(x_n)) \wedge \text{Prov}(\neg\forall x \in X\ \varphi(x))$

Subsequently, Gödel defined an ω-consistent formal system as a system that is not ω-inconsistent or a system within which both $\forall n \in \mathbb{N}\ \text{Prov}(\varphi(x_n))$ and $\text{Prov}(\neg\forall x \in X\ \varphi(x))$ cannot be true for any formula $\varphi(x)$ in the language of the formal system.

We employed a bijective function within the definition ω-inconsistency even though the function that Gödel defined to construct the Gödel numbers was an injective or one-to-one function. This is justified because a bijective function can easily be constructed from Gödel's injective function (Appendix A: Gödel Numbers and the Arithmetization of Syntax).

We will now state the customary definition of the principle of mathematical induction. For any predicate (or definite property) $P(n)$, where $n \in \mathbb{N}$:

**Principle of Mathematical Induction:** $(P(0) \wedge \forall n \in \mathbb{N}\ \{P(n) \rightarrow P(n+1)\}) \rightarrow \forall n \in \mathbb{N}\ P(n)$

This definition immediately follows from the axioms of Gödel's system, specifically axioms I.3 and IV.1 (Appendix B: Gödel's System), where the predicate $P$ is assumed to represent the formula in axiom IV.1. The principle of mathematical induction asserts that the truth of $P(0)$, called the induction basis, and the truth of $\forall n \in \mathbb{N}\ \{P(n) \rightarrow P(n+1)\}$, called the induction step, implies the truth of $\forall n \in \mathbb{N}\ P(n)$. The $P(n)$ within the induction step is called the induction hypothesis because the truth of $P(n)$ is assumed for any arbitrary $n \in \mathbb{N}$ during the proof of the induction step. Also within the induction step, the notation $n + 1$ denotes the successor operation. Although the successor operation $S$ is customarily defined as $S(x) = x + 1$, the successor operation does not simply define addition by a constant, but rather expresses the concept of the next element or successor within some complete enumeration.

Since the principle of mathematical induction only expresses the underlying concept that justifies a proof that employs the principle of mathematical induction, not its actual implementation, we will also define proof by mathematical induction within the context of Gödel's proof. To construct the definition we first substitute $\text{Prov}(\varphi(x))$ for the predicate $P(n)$ in the definition of the principle of mathematical induction, where $\varphi(x)$ denotes a formula that defines the predicate $P(n)$ within the formal system. This substitution is justified by Gödel's theorem 5 if the predicate $P$ is assumed to represent a primitive recursive relation and thus the formula $\varphi$ represents a numeric formula that implements or computes the relation. We then employ the bijective function $f$ that was previously defined to broaden the applicability of the definition, which reflects how the proof is generally applied. As a result, the following definition follows from the axioms of Gödel's formal system along with Gödel's theorem 5. For any primitive recursive numeric formula $\varphi(x)$ in some formal system and some bijective function $f$: $\mathbb{N} \rightarrow X$ and thus also $x_n = f(n)$:

**Proof by Mathematical Induction:** $(\text{Prov}(\varphi(x_0)) \wedge \forall n \in \mathbb{N}\ \{\text{Prov}(\varphi(x_n)) \rightarrow \text{Prov}(\varphi(x_{n+1}))\}) \rightarrow$
$\forall n \in \mathbb{N}\ \text{Prov}(\varphi(x_n))$

The definition asserts that if a proof of the induction basis $\varphi(x_0)$ exists and if a proof of $\varphi(x_{n+1})$ must exist whenever a proof of the induction hypothesis $\varphi(x_n)$ exists then that implies that a proof of $\varphi(x_n)$ exists for every $n \in \mathbb{N}$ or $\forall n \in \mathbb{N}\ \text{Prov}(\varphi(x_n))$, where the notation $n + 1$ again denotes the successor operation and thus $x_{n+1}$ denotes the successor of $x_n$. Within this definition the induction step $\forall n \in \mathbb{N}\ \{\text{Prov}(\varphi(x_n)) \rightarrow \text{Prov}(\varphi(x_{n+1}))\}$ consists of proving that, for any arbitrary $n \in \mathbb{N}$, a proof of the induction hypothesis $\varphi(x_n)$ can be transformed into a proof of $\varphi(x_{n+1})$.

Therefore, since the first conjunct in the definition of ω-inconsistency is the consequent in the definition of proof by mathematical induction, it is clear how the concepts of ω-inconsistency and mathematical induction may be related. Thus, at least superficially, the concept of ω-inconsistency postulates that, for some formula $\varphi(x)$ in the language of a formal system, a proof by mathematical induction that proves that the sentence $\varphi(x_n)$ is true for every $n \in \mathbb{N}$ is then seemingly contradicted by another proof that proves that the sentence $\neg\forall x \in X\ \varphi(x)$ is also true.

## 6. THE INTUITIONISTIC ACCEPTABILITY OF GÖDEL'S PROOF

Gödel stated that the proof of his first incompleteness theorem is constructive and thus also intuitionistically acceptable ([1.a], pg. 177). Therefore, in this section we will briefly examine what constitutes a constructive or intuitionistically acceptable proof. However, we will not be concerned with whether Gödel's proof can be formulated within an intuitionistic system of logic, but rather with the foundational concepts upon which Gödel based his proof and whether those concepts may be considered to be constructive.

Two primary assumptions in classical logic are considered to be intuitionistically unacceptable. In their most general form both assumptions pertain to what may be considered to be constructive and thus they also relate to what may be considered to constitute a finitary method of proof within Hilbert's program. What may be considered to be constructive has never been conclusively defined because what is constructive is inherently related to the subjective idea of what is practical or realistic, which is opposed to what is theoretical or abstract and thus also potentially impractical and unrealistic. The two intuitionistically unacceptable assumptions that will be examined are identified in Kleene ([4] pg. 46 – 53). Kleene employs the assumptions in an attempt to define a restricted form of classical logic that is intuitionistically acceptable ([4] pg. 82, 162 – 163).

The first and principal assumption of classical logic that is considered to be intuitionistically unacceptable is expressed by the formula $\neg\neg\varphi \rightarrow \varphi$, where $\varphi$ denotes any statement or sentence. The principle underlying $\neg\neg\varphi \rightarrow \varphi$ is not intuitionistically acceptable because the universal applicability of the Law of the Excluded Middle (LEM) is rejected within intuitionistic logic. The LEM asserts that, for any statement or sentence $\varphi$, either $\varphi$ or the negation of $\varphi$ is true, in symbols: $\varphi \vee \neg\varphi$. As a result, if it is not the case that the negation of $\varphi$ is true, in symbols: $\neg\neg\varphi$, then $\varphi$ must be true because $\neg\varphi$ is not true, in symbols: $\neg\neg\varphi \rightarrow \varphi$. Intuitionistic logic does not reject $\neg\neg\varphi \rightarrow \varphi$ because it is universally false, but rather because it is not necessarily true and hence possibly inapplicable depending on the context within which it is applied.

Gödel assumed $\neg\neg\varphi \rightarrow \varphi$ when he substituted the relation $Q$ for the relation in 2 of his theorem 5. The result of this substitution is $\neg\neg(\varphi_x \text{ ProofOf } \varphi_y(y)) \rightarrow \text{Prov}(\neg\varphi_q(x, y))$, due to the definition of the relation $Q$, but the result is stated as $\varphi_x \text{ ProofOf } \varphi_y(y) \rightarrow \text{Prov}(\neg\varphi_q(x, y))$ in 3 of Gödel's proof (section 4: Gödel's Proof), which is due to the assumption of $\neg\neg\varphi \rightarrow \varphi$. By assuming $\neg\neg\varphi \rightarrow \varphi$ within this context Gödel asserts that the primitive recursive relation $Q$ must be either true or false. There is no other option. As a consequence, Gödel's assumption of $\neg\neg\varphi \rightarrow \varphi$ within this context is not necessarily intuitionistically unacceptable because it is obvious that a primitive recursive relation must actually be either true or false due to how the primitive recursive functions were defined by Gödel.

The second assumption of classical logic that is considered to be intuitionistically unacceptable is related to the first and may be expressed by how the quantifiers are interpreted, defined and evaluated. Therefore, we will begin by considering two different interpretations of the existential quantifier. For any formula $\varphi(x)$, if some specific defined value has been explicitly identified for which $\varphi(x)$ is true then $\exists x\ \varphi(x)$ simply reflects that $\varphi(x)$ is true for this value. We will call this the constructive interpretation of the existential quantifier. However, if no specific defined value has been explicitly identified for which $\varphi(x)$ is true then $\exists x\ \varphi(x)$ is merely the abstract assertion that some value exists for which $\varphi(x)$ is true, although no such actual value may be known. We will call this the definitional interpretation of the existential quantifier.

The constructive interpretation of the existential quantifier requires a specific defined value to be explicitly identified for which the formula $\varphi(x)$ is true before the sentence $\exists x\ \varphi(x)$ may be asserted to be true. As a result, any proof of $\exists x\ \varphi(x)$ is required to find or construct a specific defined value for which $\varphi(x)$ can be proven to be true. In contrast, the definitional interpretation of the existential quantifier does not place any such restrictions on when the sentence $\exists x\ \varphi(x)$ may be asserted to be true nor on how it may be proven. However, in practice, if $\exists x\ \varphi(x)$ is not proven constructively with a direct proof then it is typically proven with an indirect proof or a proof that employs the principle of *reductio ad absurdum*, which then also requires $\neg\neg\varphi \rightarrow \varphi$ to be assumed to be true.

Next, we will consider the symmetrical definition of the quantifiers, which is how the quantifiers are typically defined. What we call the symmetrical definition of the quantifiers is when one of the quantifiers is taken as primitive and thus defined intuitively and then the other quantifier is formally defined in terms of it. If the existential quantifier is taken as primitive then the universal quantifier is formally defined in terms of it, in symbols: $\forall x\ \varphi(x) \equiv_{df} \neg\exists x\ \neg\varphi(x)$. Conversely, if the universal quantifier is taken as primitive then the existential quantifier is formally defined in terms of it, in symbols: $\exists x\ \varphi(x) \equiv_{df} \neg\forall x\ \neg\varphi(x)$. Gödel formally states the latter definition in his relation 32.d (Appendix B: Gödel's system). Subsequently, both of the equalities $\forall x\ \varphi(x) = \neg\exists x\ \neg\varphi(x)$ and $\exists x\ \varphi(x) = \neg\forall x\ \neg\varphi(x)$ can then be proven to be true regardless of which quantifier is taken as primitive; although this also requires $\neg\neg\varphi \rightarrow \varphi$ to be assumed to be true as well. As a result, the customary equalities, which are generically defined as $\forall = \neg\exists\neg$, $\neg\forall = \exists\neg$, $\exists = \neg\forall\neg$ and $\neg\exists = \forall\neg$, are then also true, which is why the definition is characterized as symmetrical.

The symmetrical definition of the quantifiers formally asserts the definitional interpretation of the existential quantifier to be true because the definition formally equates the two quantifiers and the universal quantifier does not have a constructive interpretation. The universal quantifier does not have a constructive interpretation because explicitly identifying an infinite number of specific defined values is clearly not practical or realistic. Although a proof by mathematical induction may be considered to be a direct proof of such a formula because a proof by mathematical induction does not necessarily employ the principle of *reductio ad absurdum*, mathematical induction is generally not considered to be constructive and thus does not justify a constructive interpretation of the universal quantifier.

Finally, we will consider how the quantifiers are evaluated or, more specifically, how any formula that contains one or more quantifiers may actually be implemented or computed. The general method that is employed to evaluate the existential or universal closure of any formula $\varphi(x)$, which is $\exists x\ \varphi(x)$ or $\forall x\ \varphi(x)$ respectively, is some type of linear search through every possible value that the variable $x$ can take. However, if the value of the variable $x$ is not explicitly bounded within the formula $\varphi(x)$ and the number of the possible values that the variable $x$ can take is not finite then employing any type of linear search to evaluate $\exists x\ \varphi(x)$ or $\forall x\ \varphi(x)$ may potentially lead to an infinite search, which is a search that never terminates. Therefore, since an infinite search signifies a result that is effectively indeterminate then the sentence $\exists x\ \varphi(x)$ or the sentence $\forall x\ \varphi(x)$ may effectively be neither true nor false even though, theoretically, every sentence must ultimately be either true or false if the LEM and thus also $\neg\neg\varphi \rightarrow \varphi$ are assumed to be true.

The intuitive distinction between the constructive and definitional interpretations of the existential quantifier may not be especially relevant within every context, however, that is certainly not the case within the context of Gödel's proof. Gödel's provability predicate is defined in relation 46 as, $\text{Prov}(\varphi_x) \equiv_{df} \exists n\ (\varphi_n \text{ ProofOf } \varphi_x)$. Subsequently, for some specific defined sentence $\varphi_m$, if $\text{Prov}(\varphi_m)$ is stated and proven, but the proof of $\text{Prov}(\varphi_m)$ does not produce an actual proof of $\varphi_m$ and does not indicate how an actual proof may be constructed then the proof of $\text{Prov}(\varphi_m)$ merely proves the abstract assertion that a proof of $\varphi_m$ must exist. However, to prove that a sentence is true it is customary to produce an actual proof of the sentence. It is not generally considered to be sufficient to prove that a proof of a sentence must exist without indicating how an actual proof can be constructed. Furthermore, an actual proof of a sentence also typically indicates or at least suggests why the sentence is actually true not merely why it must be true, which is a subtle but potentially meaningful distinction.

As a consequence, since Gödel's provability predicate clearly employs an unbounded existential quantifier then the predicate is not necessarily intuitionistically acceptable because an evaluation of the predicate may potentially result in an infinite search. Specifically, for some sentence $\varphi_m$, if $\text{Prov}(\varphi_m)$ is not true, in symbols: $\neg\exists n\ (\varphi_n \text{ ProofOf } \varphi_m)$, and the number of possible values that the variable $n$ can take is infinite then an evaluation of $\text{Prov}(\varphi_m)$ will not necessarily produce a result. Furthermore, if the symmetrical definition of the quantifiers is employed to define the quantifiers then $\neg\text{Prov}(\varphi_m) = \neg\exists n\ (\varphi_n \text{ ProofOf } \varphi_m) = \forall n\ \neg(\varphi_n \text{ ProofOf } \varphi_m)$ and thus $\neg\text{Prov}(\varphi_m)$ clearly does not have a constructive interpretation. However, Gödel is still able to assert that the proof of his first incompleteness theorem is constructive and thus also intuitionistically acceptable. This will be addressed in more detail in the next section.

## 7. A METICULOUS INTUITIVE ANALYSIS OF GÖDEL'S PROOF

We will now provide a detailed, line-by-line analysis of the proof of Gödel's first incompleteness theorem (section 4: Gödel's Proof).

### 7.1 Theorem 5

Gödel's theorem 5 states that for any primitive recursive relation $R(x_1, \ldots, x_n)$, where $n \geq 1$, a numeric formula $\varphi$ with $n$ free variables exists such that $R(x_1, \ldots, x_n) \rightarrow \text{Prov}(\varphi(x_1, \ldots, x_n))$ and $\neg R(x_1, \ldots, x_n) \rightarrow \text{Prov}(\neg\varphi(x_1, \ldots, x_n))$. The first formula states that if the relation $R$ is true for the specific numeric values $x_1, \ldots, x_n$ then a deductive proof of the numeric formula $\varphi(x_1, \ldots, x_n)$ exists that proves the truth of $R(x_1, \ldots, x_n)$ or proves that the relation $R$ is true for the specific numeric values $x_1, \ldots, x_n$. Subsequently, the second formula states the inverse, which is that if the relation $R$ is false for the specific numeric values $x_1, \ldots, x_n$ then a deductive proof of the numeric formula $\varphi(x_1, \ldots, x_n)$ exists that proves the truth of the negation of $R(x_1, \ldots, x_n)$ or proves that the relation $R$ is not true for the specific numeric values $x_1, \ldots, x_n$.

The principal role of theorem 5 within Gödel's proof, which is of crucial importance to the proof, is that the theorem proves that, for any primitive recursive relation $R$ and any specific numeric values $x_1, \ldots, x_n$, $R(x_1, \ldots, x_n)$ can be proven constructively to be either true or false. Since Gödel defined a primitive recursive relation as a relation that can be implemented or computed by a primitive recursive function then a numeric formula $\varphi$ that defines a primitive recursive function within the formal system must exist to implement or compute every primitive recursive relation. Subsequently, Gödel employed a proof by mathematical induction to prove that every primitive recursive function can be defined within his formal system. As a result, for any numeric formula $\varphi$ that implements or computes a primitive recursive relation $R$ and for any specific numeric values $x_1, \ldots, x_n$, the finite computation of $\varphi(x_1, \ldots, x_n)$, which must be finite due to how Gödel defined the primitive recursive functions, can be considered to be a constructive proof of $R(x_1, \ldots, x_n)$.

As a consequence, the finite computation of the numeric formula $\varphi$ for the specific numeric values $x_1, \ldots, x_n$ constitutes a deductive proof of $\varphi(x_1, \ldots, x_n)$ and, since $\varphi$ defines a primitive recursive function that implements or computes the primitive recursive relation $R$ in theorem 5, the finite computation of $\varphi(x_1, \ldots, x_n)$ also constitutes a deductive proof of $R(x_1, \ldots, x_n)$. Therefore, although the first formula of theorem 5 is $R(x_1, \ldots, x_n) \rightarrow \exists y\ (\varphi_y \text{ ProofOf } \varphi(x_1, \ldots, x_n))$ and the second formula is $\neg R(x_1, \ldots, x_n) \rightarrow \exists y\ (\varphi_y \text{ ProofOf } \neg\varphi(x_1, \ldots, x_n))$, due to definition of Gödel's provability predicate, the unbounded existential quantifier in the formulas will not result in an infinite search because the proof of theorem 5 along with the definition of the primitive recursive functions indicates precisely how an actual deductive proof of $\varphi(x_1, \ldots, x_n)$ and thus also $R(x_1, \ldots, x_n)$ can be constructed. Hence, a search is not required and thus an infinite search will never result.

### 7.2 The Relation *Q*

To begin his proof Gödel defined the relation $Q$ as $Q(x, y) \equiv_{\text{df}} \neg(\varphi_x \text{ ProofOf } \varphi_y(y))$, which states that the relation $Q(x, y)$ is true for the values $x$ and $y$ if it is not the case that $x = [\varphi_x]$ is the Gödel number of a finite sequence of formulas and $\varphi_x$ constitutes a deductive proof of the sentence $\varphi_y(y)$. The sentence $\varphi_y(y)$ is the result of substituting the value $y = [\varphi_y(x)]$ for the single free variable in the formula $\varphi_y(x)$. Hence, the Gödel number of the formula $\varphi_y(x)$ serves as a fixed value, often called a fixed point, with which to evaluate the formula $\varphi_y(x)$.

A formula with a single free variable is what Gödel called a class expression. A formula with a single free variable, such as $\varphi_y(x)$, cannot be determined to be either true or false until some value has been substituted for its free variable. For example, the formula "$x + 2 = 3$" is neither true nor false until some specific number is substituted for the variable $x$. Thus, if $x = 1$ then "$1 + 2 = 3$" is true and if $x = 2$ then "$2 + 2 = 3$" is false. Thus, Gödel called a formula with a single free variable a class expression because it segregates every possible value of its single free variable into two distinct classes or sets of values; the set that contains every value for which the formula is true and the set that contains every value for which the formula is false.

Hence, Gödel chose to use the Gödel number of the formula $\varphi_y(x)$ itself to serve as a fixed value for which the formula could be determined to be either true or false. This may seem arbitrary or even contrived because then only the sentences $\varphi_0(0)$, $\varphi_1(1)$, $\varphi_2(2)$, … would be evaluated by the relation $Q$. But, Gödel could have defined the relation $Q$ with three free variables instead of two, in symbols: $Q(x, y, z) \equiv_{\text{df}} \neg(\varphi_x \text{ ProofOf } \varphi_y(z))$, which states that the relation $Q(x, y, z)$ is true if it is not the case that $x = [\varphi_x]$ is the Gödel number of a finite sequence of formulas and $\varphi_x$ constitutes a deductive proof of the formula $y = [\varphi_y(x)]$ when the value $z$ is substituted for the single free variable in the formula $\varphi_y(x)$. However, since Gödel only needed to evaluate the relation $Q(x, y, z)$ when $z = y$ to carry out his proof, which would be $Q(x, y, y)$, he simply defined the relation $Q$ with two free variables.

Since the relation ProofOf is primitive recursive and thus the relation $Q$ is also primitive recursive, Gödel then substitutes the primitive recursive relation $Q$ for the relation in theorem 5. As a result, there exists a numeric formula that implements or computes the primitive recursive relation $Q$, which is denoted by $\varphi_q$, and thus:

**1:** $Q(x, y) = \neg(\varphi_x \text{ ProofOf } \varphi_y(y)) \rightarrow \text{Prov}(\varphi_q(x, y))$

**2:** $\neg Q(x, y) = \neg\neg(\varphi_x \text{ ProofOf } \varphi_y(y)) = \varphi_x \text{ ProofOf } \varphi_y(y) \rightarrow \text{Prov}(\neg\varphi_q(x, y))$

$Q(x, y) = \neg(\varphi_x \text{ ProofOf } \varphi_y(y))$ due to the definition of the relation $Q$. Hence, 1 states that if the relation $Q$ is true for the numeric values $x$ and $y$ then a proof of the numeric formula $\varphi_q(x, y)$ exists that proves that $\varphi_q(x, y)$ is true for those values. Similarly, $\neg Q(x, y) = \varphi_x \text{ ProofOf } \varphi_y(y)$ due to the definition of the relation $Q$ along with $\neg\neg\varphi \rightarrow \varphi$. Hence, 2 states that if the relation $Q$ is not true for the numeric values $x$ and $y$ then a proof of the numeric formula $\varphi_q(x, y)$ exists that proves that $\varphi_q(x, y)$ is not true for those values. The proof of the numeric formula $\varphi_q(x, y)$ in both 1 and 2 is the finite computation of $\varphi_q(x, y)$ that constitutes a deductive proof of $\varphi_q(x, y)$, which must exist because the numeric formula $\varphi_q(x, y)$ defines a primitive recursive function.

### 7.3 The Relations *P* and *R*

Gödel defined the relation $P$ in terms of the relation $Q$ as $\varphi_p(x) \equiv_{\text{df}} \forall n\ \varphi_q(n, x)$. Thus, if the variable $x$ is replaced by $y$ for clarity then $P(y) \equiv_{\text{df}} \forall n\ \neg(\varphi_n \text{ ProofOf } \varphi_y(y))$ due to the definition of the relation $Q$. Hence, the definition states that the relation $P(y)$ is true for the value $y$ if, for every possible value of the variable $n$, it is not the case that both $n = [\varphi_n]$ is the Gödel number of a finite sequence of formulas and $\varphi_n$ constitutes a deductive proof of the sentence $\varphi_y(y)$. Again, the sentence $\varphi_y(y)$ is the result of substituting the value $y = [\varphi_y(x)]$ for the single free variable in the formula $\varphi_y(x)$. The relation intuitively states that $P(y)$ is true for any numeric value $y$ if a deductive proof of the sentence $\varphi_y(y)$ does not exist. We also note that, like Gödel's provability predicate, the relation $P$ is not primitive recursive.

Gödel then defined the relation $R$ in terms of both of the relations $Q$ and $P$ as $\varphi_r(x) \equiv_{\text{df}} \varphi_q(x, p)$, where $p = [\varphi_p(x)]$. Thus, $R(x) \equiv_{\text{df}} \neg(\varphi_x \text{ ProofOf } \varphi_p(p))$ due to the definition of the relation $Q$. Hence, the definition states that the relation $R(x)$ is true for the value $x$ if it is not the case that both $x = [\varphi_x]$ is the Gödel number of a finite sequence of formulas and $\varphi_x$ constitutes a deductive proof of the sentence $\varphi_p(p)$. The sentence $\varphi_p(p)$ is the result of substituting the value $p = [\varphi_p(x)]$ for the single free variable in the formula $\varphi_p(x)$. As a result, $\varphi_p(p)$ denotes the relation $P$ with its own Gödel number substituted for its single free variable. The relation intuitively states that $R(x)$ is true for any numeric value $x$ if $x = [\varphi_x]$ does not constitute a deductive proof of the sentence $\varphi_p(p)$.

Gödel then states the equalities $\varphi_q(x, p) = \varphi_r(x)$ and $\varphi_p(p) = \forall n\ \varphi_q(n, p) = \forall n\ \varphi_r(n)$. The equalities immediately follow from the definitions of the relations $P$ and $R$. The equality $\varphi_q(x, p) = \varphi_r(x)$ states that the relation $Q(x, p)$, which is the relation $Q(x, y)$ when the value of its second free variable is $p = [\varphi_p(x)]$, is identical to the relation $R(x)$. Thus, $\varphi_q(x, p)$ is true for the value $x$ if and only if $\varphi_r(x)$ is true for the value $x$. The equality $\varphi_p(p) = \forall n\ \varphi_q(n, p) = \forall n\ \varphi_r(n)$ states that the relation $P(p)$, which is the relation $P(x)$ when the value of its single free variable is $p = [\varphi_p(x)]$, is identical to the universal closure of the relation $R(x)$ or $\forall n\ R(n)$. Thus, $\varphi_p(p)$ is true if and only if $\forall n\ \varphi_r(n)$ is true.

### 7.4 The Result

Gödel then derives the two formulas upon which the result of the proof of his first incompleteness theorem is based. He does so by substituting the value $p = [\varphi_p(x)]$ for the free variable $y$ in 1 and 2 (section 7.2: The Relation $Q$), which are 2 and 3 in his proof (section 4: Gödel's Proof). The result of this substitution is:

**3:** $\neg(\varphi_x \text{ ProofOf } \varphi_p(p)) \rightarrow \text{Prov}(\varphi_q(x, p))$

**4:** $\varphi_x \text{ ProofOf } \varphi_p(p) \rightarrow \text{Prov}(\neg\varphi_q(x, p))$

And then applying the equality $\varphi_p(p) = \forall n\ \varphi_r(n)$ to the relation in the antecedent and applying the equality $\varphi_q(x, p) = \varphi_r(x)$ to the relation within the provability predicate in the consequent:

**5:** $\neg(\varphi_x$ ProofOf $\forall n\ \varphi_r(n)) \rightarrow$ Prov($\varphi_r(x)$)

**6:** $\varphi_x$ ProofOf $\forall n\ \varphi_r(n) \rightarrow$ Prov($\neg\varphi_r(x)$)

As a result, 5 and 6 are 8 and 9 in Gödel's proof (section 4: Gödel's Proof). 5 states that if $\varphi_x$ is not a deductive proof of $\forall n\ \varphi_r(n)$ then a deductive proof of the numeric formula $\varphi_r(x)$ exists for the value $x$, which comports with the definition of the relation because $\varphi_r(x)$ is true for any numeric value $x$ if $x = [\varphi_x]$ does not constitute a deductive proof of the sentence $\varphi_p(p) = \forall n\ \varphi_r(n)$. And 6 states that if $\varphi_x$ is a deductive proof of $\forall n\ \varphi_r(n)$ then a deductive proof of numeric formula $\neg\varphi_r(x)$ exists for the value $x$, which also comports with the definition of the relation because $\varphi_r(x)$ is false for any numeric value $x$ if $x = [\varphi_x]$ constitutes a deductive proof of the sentence $\varphi_p(p) = \forall n\ \varphi_r(n)$.

Gödel then describes the contradictions that result from 5 and 6 intuitively. However, we will formally derive the contradictions that Gödel intuitively inferred from 5 and 6. First, assume some $x = [\varphi_x]$ exists such that $\varphi_x$ is a deductive proof of $\forall n\ \varphi_r(n)$:

**7:** $\exists x\ (\varphi_x$ ProofOf $\forall n\ \varphi_r(n))$

Due to the definition of Gödel's provability predicate we can immediately infer from 7:

**8:** Prov($\forall n\ \varphi_r(n)$)

And then from 6 and 7, using modus ponens:

**9:** $\exists x$ Prov($\neg\varphi_r(x)$)

But, if the formal system is simply consistent then 9 contradicts 8. Therefore, assume a deductive proof of $\neg\forall n\ \varphi_r(n)$ exists and thus, due to the definition of Gödel's provability predicate:

**10:** Prov($\neg\forall n\ \varphi_r(n)$)

Since the formal system is assumed to be simply consistent then it is also assumed that a $x = [\varphi_x]$ does not exist such that $\varphi_x$ is a deductive proof of $\forall n\ \varphi_r(n)$, which is the refutation of $\neg\forall n\ \varphi_r(n)$:

**11:** $\forall x\ \neg(\varphi_x$ ProofOf $\forall n\ \varphi_r(n))$

And then from 5 and 11, again using modus ponens:

**12:** $\forall x$ Prov($\varphi_r(x)$)

But, if the formal system is ω-consistent then 12 contradicts 10. ■

It is first assumed that some $x = [\varphi_x]$ exists such that $\varphi_x$ is a deductive proof of $\forall n\ \varphi_r(n)$, which is stated in 7. But this assumption results in a contradiction if the formal system is simply consistent. Therefore, it is assumed that a deductive proof of $\neg\forall n\ \varphi_r(n)$ exists, which is stated in 10, and since the formal system is assumed to be simply consistent and $\forall n\ \varphi_r(n)$ is the refutation of $\neg\forall n\ \varphi_r(n)$ then it is also assumed that a $x = [\varphi_x]$ does not exist such that $\varphi_x$ is a deductive proof of $\forall n\ \varphi_r(n)$, which is stated in 11. But the assumption that a deductive proof of $\neg\forall n\ \varphi_r(n)$ exists and the further assumption that a deductive proof of its refutation, $\forall n\ \varphi_r(n)$, does not exist also results in a contradiction if the formal system is ω-consistent.

As a result, Gödel infers that $\forall n\ \varphi_r(n)$ cannot be proven to be either true or false within the formal system because neither a deductive proof of $\forall n\ \varphi_r(n)$ nor a deductive proof of its refutation can exist within the formal system if the formal system is simply consistent as well as ω-consistent. Gödel calls $\forall n\ \varphi_r(n)$ an undecidable sentence because the truth of $\forall n\ \varphi_r(n)$ cannot be decided within the formal system and thus he also calls the formal system incomplete because the system contains an undecidable sentence. Hence, Gödel further asserts that any formal system that is robust enough to include elementary arithmetic and mathematical logic must be inherently incomplete because the foundation of the formal system that Gödel employed to derive his result is constructed using only the most basic and widely accepted axioms and methods of elementary arithmetic and mathematical logic.

Gödel further states that $\forall x\ \varphi_r(x)$ can easily be seen to be true precisely because $\varphi_p(p) = \forall n\ \varphi_r(n)$ is an undecidable sentence and thus cannot be proven within the formal system. Since $\forall n\ \varphi_r(n)$ cannot be proven within the formal system then a deductive proof of neither $\forall n\ \varphi_r(n)$ or its refutation exists within the formal system and thus $\forall x\ \neg(\varphi_x \text{ ProofOf } \forall n\ \varphi_r(n))$ is clearly true. However, $\forall x\ \varphi_r(x) = \forall x\ \neg(\varphi_x \text{ ProofOf } \forall n\ \varphi_r(n))$, due to the definition of the formula $\varphi_r$, the definition of the relation $Q$ and the equality $\varphi_p(p) = \forall n\ \varphi_r(n)$. Hence, $\forall x\ \varphi_r(x)$ is clearly true.

Since a formal system is intuitively defined to be complete if every meaningful sentence in the language of the formal system can be proven within the formal system, the importance of the assertion that $\forall x\ \varphi_r(x)$ is actually a true sentence is that the assertion refutes the potential criticism that the sentence $\forall x\ \varphi_r(x)$ expresses a contrived or nonsensical statement and is thus meaningless.

### 7.5 The Contradictions

We will begin by defining $\bot_1 \equiv_{\text{df}} \exists x\ (\varphi_x \text{ ProofOf } \forall n\ \varphi_r(n)) \rightarrow \text{Prov}(\forall n\ \varphi_r(n)) \wedge \exists x\ \text{Prov}(\neg\varphi_r(x))$ and $\bot_2 \equiv_{\text{df}} \text{Prov}(\neg\forall n\ \varphi_r(n)) \rightarrow \forall x\ \neg(\varphi_x \text{ ProofOf } \forall n\ \varphi_r(n)) \wedge \forall x\ \text{Prov}(\varphi_r(x))$. $\bot_1$ denotes the first contradiction, which is the contradiction that is expressed in 7, 8 and 9, and $\bot_2$ denotes the second contradiction, which is the contradiction that is expressed in 10, 11 and 12.

Gödel first interprets $\text{Prov}(\forall n\ \varphi_r(n)) \wedge \exists x\ \text{Prov}(\neg\varphi_r(x))$ in $\bot_1$ as a contradiction if the formal system is simply consistent. This is justified because $\exists x\ \text{Prov}(\neg\varphi_r(x))$ states that a proof of $\neg\varphi_r(x)$ exists in the theory of the formal system for some fixed value that is a valid substitution for its free variable and thus, for some fixed value $c$, $\exists x\ \text{Prov}(\neg\varphi_r(x)) \rightarrow \text{Prov}(\neg\varphi_r(c))$ may be inferred. As a consequence, since $\text{Prov}(\forall n\ \varphi_r(n)) \rightarrow \text{Prov}(\varphi_r(c))$ for any fixed value $c$, due to axiom III.1 (Appendix B: Gödel's System), then a deductive proof exists within the formal system for both of the sentences $\neg\varphi_r(c)$ and $\varphi_r(c)$ for some fixed value $c$ and thus the formal system is simply inconsistent if $\text{Prov}(\forall n\ \varphi_r(n)) \wedge \exists x\ \text{Prov}(\neg\varphi_r(x))$ is true.

Gödel then interprets $\text{Prov}(\neg\forall n\ \varphi_r(n))$ and $\forall x\ \text{Prov}(\varphi_r(x))$ in $\bot_2$ as a contradiction if the formal system is $\omega$-consistent. This is justified because Gödel specifically defined $\text{Prov}(\neg\forall n\ \varphi_r(n))$ and $\forall x\ \text{Prov}(\varphi_r(x))$ to be an $\omega$-inconsistency. Thus, the pertinent question is not how the contradiction is inferred, but why the assumption that the formal system is simply consistent is not sufficient to interpret $\text{Prov}(\neg\forall n\ \varphi_r(n))$ and $\forall x\ \text{Prov}(\varphi_r(x))$ as a contradiction.

If $\neg\neg\varphi \rightarrow \varphi$ is assumed to be true then $\text{Prov}(\neg\forall n\ \varphi_r(n)) = \text{Prov}(\exists n\ \neg\varphi_r(n))$, due to relation 32.d (Appendix B: Gödel's System). Hence, if the existential quantifier is interpreted constructively then some fixed value $c$ must exist such that $\text{Prov}(\neg\forall n\ \varphi_r(n)) \rightarrow \text{Prov}(\neg\varphi_r(c))$. Subsequently, since $\forall x\ \text{Prov}(\varphi_r(x))$ states that $\varphi_r(x)$ exists in the theory of the formal system for any fixed value that is a valid substitution for its free variable then, for any fixed value $c$, $\forall x\ \text{Prov}(\varphi_r(x)) \rightarrow \text{Prov}(\varphi_r(c))$ may be inferred. As a consequence, a deductive proof again exists within the formal system for both of the sentences $\neg\varphi_r(c)$ and $\varphi_r(c)$ for some fixed value $c$ and thus the formal system is simply inconsistent if both $\text{Prov}(\neg\forall n\ \varphi_r(n))$ and $\forall x\ \text{Prov}(\varphi_r(x))$ are true. Therefore, if the assumption that the formal system is simply consistent is not sufficient to infer a contradiction from the conjunction of $\text{Prov}(\neg\forall n\ \varphi_r(n))$ and $\forall x\ \text{Prov}(\varphi_r(x))$ then either $\text{Prov}(\neg\forall n\ \varphi_r(n)) \rightarrow \text{Prov}(\neg\varphi_r(c))$ or $\forall x\ \text{Prov}(\varphi_r(x)) \rightarrow \text{Prov}(\varphi_r(c))$ or both cannot be inferred from the axioms of the formal system for any fixed value $c$.

## 8. ROSSER'S PROOF

We will first define Rosser's proof relation and provability predicate.

**Rosser's Proof Relation:** $\varphi_x \text{ Proves } \varphi_y \equiv_{\text{df}} \varphi_x \text{ ProofOf } \varphi_y \wedge \neg\exists z\ (z \leq x \wedge \varphi_z \text{ ProofOf } \neg\varphi_y)$

The validity of Rosser's proof relation immediately follows from the assumption that the formal system is simply consistent along with the validity of Gödel's proof relation, which is employed within the definition. The definition states that the relation $\varphi_x$ Proves $\varphi_y$ is true if $\varphi_x$ constitutes a deductive proof of $\varphi_y$ and a $z = [\varphi_z]$ does not exist such that $z$ is smaller than or equal to $x = [\varphi_x]$ and $\varphi_z$ constitutes a deductive proof of the negation of $\varphi_y$. The value of the variable $z$ is explicitly bounded within the definition so that, like Gödel's proof relation, Rosser's proof relation is also primitive recursive.

**Rosser's Provability Predicate:** $\text{Prov}(\varphi_x) \equiv_{\text{df}} \exists n\ (\varphi_n \text{ Proves } \varphi_x)$

The definition of Rosser's provability predicate is essentially identical to the definition of Gödel's provability predicate with the obvious exception that Rosser's proof relation is employed within the definition instead of Gödel's proof relation. Thus, like Gödel's provability predicate, Rosser's provability predicate is also not primitive recursive.

The following formalization of Rosser's proof is primarily due to Kleene ([4] pg. 208). We use Kleene's formalization of Rosser's proof because Rosser did not actually construct a formal proof, but merely implied how one could be constructed, whereas Kleene constructed a detailed and complete formal proof. An asterisk followed by a number refers to a formula that Kleene employs as an assumption within Rosser's proof. These formulas will be listed without proof following Rosser's proof.

We will begin Rosser's proof by defining the relation $Q$. The definition of the relation is essentially identical to how the relation is defined in Gödel's proof except that Rosser's proof relation is employed within the definition instead of Gödel's proof relation.

**1:** $Q(x, y) \equiv_{\mathrm{df}} \neg(\varphi_x \text{ ProofOf } \varphi_y(y) \wedge \neg\exists z\ \{z \leq x \wedge \varphi_z \text{ ProofOf } \neg\varphi_y(y)\})$

Again, due to Gödel's theorem 5, there is a numeric formula $\varphi_q(x, y)$ that implements or computes the relation $Q$. The formula $\varphi_p$ is then defined identically to how it was defined in Gödel's proof:

**2:** $\varphi_p(x) \equiv_{\mathrm{df}} \forall n\ \varphi_q(n, x)$

Substituting $p = [\varphi_p(x)]$ for $x$ in 2 and then due to the definition of the relation $Q$ along with *57b and $\neg\neg\varphi \rightarrow \varphi$:

**3:** $\varphi_p(p) = \forall n\ \varphi_q(n, p) = \forall n\ (\neg\{\varphi_n \text{ ProofOf } \varphi_p(p)\} \vee \exists z\ \{z \leq n \wedge \varphi_z \text{ ProofOf } \neg\varphi_p(p)\})$

Now, assume $\varphi_p(p)$ and thus also $\text{Prov}(\varphi_p(p))$. Hence, there exists some fixed $k = [\varphi_k]$ such that:

**4:** $\varphi_k \text{ ProofOf } \varphi_p(p)$

And due to Rosser's proof relation (and the assumption of simple consistency):

**5:** $\neg(\varphi_0 \text{ ProofOf } \neg\varphi_p(p)), \neg(\varphi_1 \text{ ProofOf } \neg\varphi_p(p)), \ldots, \neg(\varphi_k \text{ ProofOf } \neg\varphi_p(p))$

Then, due to *166a:

**6:** $\forall z\ (z \leq k \rightarrow \neg(\varphi_z \text{ ProofOf } \neg\varphi_p(p))$

Due to 4 and 6, using $\exists$ and $\wedge$ introduction:

**7:** $\exists x\ (\varphi_x \text{ ProofOf } \varphi_p(p) \wedge \forall z\ \{z \leq x \rightarrow \neg(\varphi_z \text{ ProofOf } \neg\varphi_p(p))\})$

Due to *58b and *86:

**8:** $\exists x\ (\varphi_x \text{ ProofOf } \varphi_p(p) \wedge \neg\exists z\ \{z \leq x \wedge \varphi_z \text{ ProofOf } \neg\varphi_p(p)\})$

Due to *57b:

**9:** $\exists x\ \neg(\neg\{\varphi_x \text{ ProofOf } \varphi_p(p)\} \vee \exists z\ \{z \leq x \wedge \varphi_z \text{ ProofOf } \neg\varphi_p(p)\})$

Finally, due to *85a:

**10:** $\neg\forall x\ (\neg\{\varphi_x \text{ ProofOf } \varphi_p(p)\} \vee \exists z\ \{z \leq x \wedge \varphi_z \text{ ProofOf } \neg\varphi_p(p)\})$

But, due to 3, this contradicts the assumption of $\varphi_p(p)$. Therefore, assume $\neg\varphi_p(p)$ and thus also $\text{Prov}(\neg\varphi_p(p))$. Hence, there exists some fixed $k = [\varphi_k]$ such that:

**11:** $\varphi_k \text{ ProofOf } \neg\varphi_p(p)$

Then, due to *168:

**12:** $\forall x\ (x \geq k \rightarrow \exists z\ \{z \leq x \wedge \varphi_z \text{ ProofOf } \neg\varphi_p(p)\})$

And again due to Rosser's proof relation (and the assumption of simple consistency):

**13:** $\neg(\varphi_0 \text{ ProofOf } \varphi_p(p)), \neg(\varphi_1 \text{ ProofOf } \varphi_p(p)), \ldots, \neg(\varphi_{k-1} \text{ ProofOf } \varphi_p(p))$

Then, due to *166:

**14:** $\forall x\ (x < k \rightarrow \neg\{\varphi_x \text{ ProofOf } \varphi_p(p)\})$

Due to 12 and 14 and also *169:

**15:** $\forall x\ (\neg\{\varphi_x \text{ ProofOf } \varphi_p(p)\} \vee \exists z\ \{z \leq x \wedge \varphi_z \text{ ProofOf } \neg\varphi_p(p)\})$

But, due to 3, this contradicts the assumption of $\neg\varphi_p(p)$. ■

Due to 10 and 15 it is inferred that $\varphi_p(p)$ cannot be proven within the formal system and thus the formal system is incomplete because $\varphi_p(p)$ is an undecidable sentence within the formal system. Hence, the customary interpretation of Rosser's proof asserts that the result of Rosser's proof is identical to the result of Gödel's proof of his first incompleteness theorem. Subsequently, since Rosser only required the assumption that the system is simply consistent then it is further asserted that Rosser's proof renders the assumption in Gödel's proof that the formal system is ω-consistent superfluous to a proof of incompleteness.

We will now provide the list of formulas that were referenced within the preceding proof. The proofs of these formulas are due to Kleene ([4] pg. 118 – 198). Although we note that the proofs use Kleene's system, not Gödel's. The notation $\vdash$ denotes the proof-theoretic relation of deductive proof or, for some set of formulas $\Sigma$ and some single formula $\varphi$, the notation $\Sigma \vdash \varphi$ states that the formula $\varphi$ can be deductively proven from the formula or formulas in $\Sigma$.

For any formulas $\varphi$ and $\Psi$, where $x$ and $y$ are distinct variables, $k$ is some fixed numeric value, and $t$ is a term that does not contain the variable $x$:

*57b. $\varphi \wedge \neg\Psi \rightarrow \neg(\neg\varphi \vee \Psi)$
*58b. $\varphi \rightarrow \neg\Psi \leftrightarrow \neg(\varphi \wedge \Psi)$
*85a. $\exists x\, \neg\varphi(x) \rightarrow \neg\forall x\, \varphi(x)$
*86. $\neg\exists x\, \varphi(x) \leftrightarrow \forall x\, \neg\varphi(x)$
*166. $\varphi(0), \varphi(1), \ldots, \varphi(k-1) \vdash \forall x\, (x < k \rightarrow \varphi(x))$
*166a. $\varphi(0), \varphi(1), \ldots, \varphi(k) \vdash \forall x\, (x \leq k \rightarrow \varphi(x))$
*168. $\varphi(t) \vdash \forall x\, (x \geq t \rightarrow \exists y\, [y \leq x \wedge \varphi(y))$
*169. $\forall x\, (x < t \rightarrow \varphi(x)), \forall x\, (x \geq t \rightarrow \Psi(x)) \vdash \forall x\, (\varphi(x) \vee \Psi(x))$

## 9. THE DIFFERENCE BETWEEN GÖDEL'S PROOF AND ROSSER'S PROOF

The principal difference between Gödel's proof and Rosser's proof is customarily attributed to the difference between the Gödel sentence $\varphi_p(p) = \forall n\, \neg(\varphi_n \text{ ProofOf } \varphi_p(p))$ that Gödel employed in his proof and the Rosser sentence $\varphi_p(p) = \forall n\, \neg(\varphi_n \text{ ProofOf } \varphi_p(p) \wedge \neg\exists z\, \{z \leq n \wedge \varphi_z \text{ ProofOf } \neg\varphi_p(p)\})$ that Rosser employed in his proof. As a result, it is assumed that the Rosser sentence is responsible for Rosser being able to prove the incompleteness of the formal system with only the assumption that the system is simply consistent, which apparently dispenses with the assumption in Gödel's proof that the formal system is ω-consistent. The basis for this interpretation is due to Rosser's assertion that his proof relation possesses properties that Gödel's proof relation does not.

However, we will show that the crucial difference between Gödel's proof and Rosser's proof cannot be attributed to the difference between the Gödel sentence and the Rosser sentence. This will be accomplished by proving Gödel's formal system to be incomplete while employing only the Gödel sentence along with an additional explicit assumption that the formal system is simply consistent. As a result, the difference between the proofs must be attributed to other assumptions contained in Rosser's proof.

If a formal system is simply consistent then, for any sentence $\varphi$, both a proof of $\varphi$ and a proof of its negation $\neg\varphi$ cannot exist. Gödel ([1] pg. 61) and Rosser ([3.b] pg. 231) both state this definition of simple consistency. Hence, assuming $\neg(\varphi \wedge \Psi) = \neg\varphi \vee \neg\Psi$, the simple consistency of a formal system may be defined generally as $\neg(\text{Prov}(\varphi) \wedge \text{Prov}(\neg\varphi)) = \neg\text{Prov}(\varphi) \vee \neg\text{Prov}(\neg\varphi)$. Therefore, due to the definition of Gödel's provability predicate:

**Simple Consistency:** $\neg\exists x\, (\varphi_x \text{ ProofOf } \varphi) \vee \neg\exists x\, (\varphi_x \text{ ProofOf } \neg\varphi)$

Subsequently, by employing this definition of simple consistency:

**Theorem 1:** Gödel's assumption that the formal system is ω-consistent can be directly replaced by the explicit assumption of the simple consistency of the system.

*Proof*: we will employ only Gödel's formal system as it was defined in his proof along with the formal definition of simple consistency that was just defined.

We will begin by defining the relation $Q$ and the formula $\varphi_p$ identically to how they were defined in Gödel's proof (section 4: Gödel's Proof):

**1:** $Q(x, y) \equiv_{\text{df}} \neg(\varphi_x \text{ ProofOf } \varphi_y(y))$

Again, there is a numeric formula $\varphi_q(x, y)$ that implements or computes the relation $Q$:

**2:** $\varphi_p(x) \equiv_{\mathrm{df}} \forall n\ \varphi_q(n, x)$

Substituting $p = [\varphi_p(x)]$ for $x$ in 2 and then due to the definition of the relation $Q$:

**3:** $\varphi_p(p) = \forall n\ \varphi_q(n, p) = \forall n\ \neg(\varphi_n \text{ ProofOf } \varphi_p(p))$

Assume $\varphi_p(p)$ and thus also $\mathrm{Prov}(\varphi_p(p))$. Hence, due to the definition of the provability predicate:

**4:** $\exists x\ (\varphi_x \text{ ProofOf } \varphi_p(p))$

Due to relation 32.d (Appendix B: Gödel's System):

**5:** $\neg\forall x\ \neg(\varphi_x \text{ ProofOf } \varphi_p(p))$

But, due to 3, this contradicts the assumption of $\varphi_p(p)$. Therefore, assume $\neg\varphi_p(p)$ and thus also $\mathrm{Prov}(\neg\varphi_p(p))$. Hence, due to the definition of the provability predicate:

**6:** $\exists x\ (\varphi_x \text{ ProofOf } \neg\varphi_p(p))$

Thus, due to 6 and the definition of simple consistency:

**7:** $\neg\exists x\ (\varphi_x \text{ ProofOf } \varphi_p(p))$

And again due to relation 32.d:

**8:** $\neg\neg\forall x\ \neg(\varphi_x \text{ ProofOf } \varphi_p(p))$

And then assuming $\neg\neg\varphi \rightarrow \varphi$:

**9:** $\forall x\ \neg(\varphi_x \text{ ProofOf } \varphi_p(p))$

But, due to 3, this contradicts the assumption of $\neg\varphi_p(p)$. ■

The proof can be intuitively stated simply and directly. The sentence $\varphi_p(p)$ is equivalent to the sentence $\forall x\ \neg(\varphi_x \text{ ProofOf } \varphi_p(p))$ by definition, due to 3. Hence, the sentence $\varphi_p(p)$ states that a proof of itself does not exist. Consequently, if it is assumed that $\varphi_p(p)$ is true and thus also that a proof of $\varphi_p(p)$ exists then $\varphi_p(p)$ is false by definition, which is a contradiction. Subsequently, if it is assumed that $\neg\varphi_p(p)$ is true and thus also that a proof of $\neg\varphi_p(p)$ exists then a proof of $\varphi_p(p)$ cannot exist, due to the simple consistency of the system, but if a proof of $\varphi_p(p)$ does not exist then $\varphi_p(p)$ is true by definition, which is again a contradiction.

We will now examine Rosser's actual proof ([3.b] pg. 233 – 234) to further substantiate our assertion that the difference between Gödel's proof and Rosser's proof cannot be attributed to the difference between the Gödel sentence and the Rosser sentence. Rosser stated in his proof that if his Rosser sentence is substituted for the Gödel sentence in Gödel's proof then Gödel's proof can be carried out exactly as Gödel did to construct an undecidable sentence within the formal system. However, if this is actually done and $(\varphi \wedge \neg\Psi) \rightarrow \neg(\neg\varphi \vee \Psi)$ along with $\neg\neg\varphi \rightarrow \varphi$ are assumed to be true then 8 and 9 in Gödel's proof (section 4: Gödel's Proof) become:

**8:** $(\neg\{\varphi_x \text{ ProofOf } \forall n\ \varphi_r(n)\} \vee \exists z\ \{z \leq x \wedge \varphi_z \text{ ProofOf } \neg\forall n\ \varphi_r(n)\}) \rightarrow \mathrm{Prov}(\varphi_r(x))$

**9:** $(\varphi_x \text{ ProofOf } \forall n\ \varphi_r(n) \wedge \neg\exists z\ \{z \leq x \wedge \varphi_z \text{ ProofOf } \neg\forall n\ \varphi_r(n)\}) \rightarrow \mathrm{Prov}(\neg\varphi_r(x))$

It is first assumed that a deductive proof of $\forall n\ \varphi_r(n)$ exists within the formal system and also that a deductive proof of $\neg\forall n\ \varphi_r(n)$ does not exist within the formal system with a Gödel number less than or equal to the Gödel number of the deductive proof of $\forall n\ \varphi_r(n)$. As a result, the antecedent of 9 can clearly be seen to be true for some $x = [\varphi_x]$ and thus $\mathrm{Prov}(\neg\varphi_r(x))$ must be true for some value of the variable $x$. Subsequently, if the formal system is assumed to be simply consistent then $\mathrm{Prov}(\forall n\ \varphi_r(n)) \wedge \exists x\ \mathrm{Prov}(\neg\varphi_r(x))$ can be interpreted as a contradiction in precisely the same manner that it was interpreted as a contradiction in Gödel's proof.

It is then assumed that a deductive proof of $\neg\forall n\ \varphi_r(n)$ exists within the formal system and also that a deductive proof of $\forall n\ \varphi_r(n)$ does not exist within the formal system with a Gödel number less than the Gödel number of the deductive proof of $\neg\forall n\ \varphi_r(n)$. Hence, $\varphi_z \text{ ProofOf } \neg\forall n\ \varphi_r(n)$ is true for some $z = [\varphi_z]$ and thus $z \leq x \wedge \varphi_z \text{ ProofOf } \neg\forall n\ \varphi_r(n)$ is true for every value of the variable $x$ that is greater than or equal to this $z$ and $\neg(\varphi_x \text{ ProofOf } \forall n\ \varphi_r(n))$ is true for every value of the variable $x$ that is smaller than this $z$. As a result, the antecedent of 8 can clearly be seen to be true for every value of the variable $x$ and thus $\text{Prov}(\varphi_r(x))$ must also be true for every value of the variable $x$. However, $\text{Prov}(\neg\forall n\ \varphi_r(n)) \wedge \forall x\ \text{Prov}(\varphi_r(x))$ is still an $\omega$-inconsistency and thus, unless something more is assumed, the assumption that the formal system is $\omega$-consistent is still required to interpret $\text{Prov}(\neg\forall n\ \varphi_r(n)) \wedge \forall x\ \text{Prov}(\varphi_r(x))$ as contradiction.

The crucial difference between Gödel's proof and Rosser's proof that allowed Rosser to dispense with the assumption in Gödel's proof that the formal system is $\omega$-consistent is Rosser's concept of what constitutes a proof. Rosser employed a metamathematical concept of validity in his proof that necessarily supplanted the concept of proof-theoretic deductive proofs that was employed by Gödel in his proof. The concept of validity that Rosser employed in his proof may be clearly illustrated by examining the properties of Rosser's proof relation.

In Rosser's proof he stated that his proof relation possesses certain properties that Gödel's proof relation does not posses and that these properties allow him to dispense with Gödel's assumption that the formal system is $\omega$-consistent. The specific properties were formally defined by Rosser as $\text{Prov}(\varphi_n) \rightarrow \text{Prov}(\varphi_m)$ and $\text{Prov}(\neg\varphi_n) \rightarrow \text{Prov}(\neg\varphi_m)$, where $\varphi_n$ is any sentence and $m = [\text{Prov}(\varphi_n)]$. Thus, $\varphi_m$ denotes the sentence $\text{Prov}(\varphi_n)$. Subsequently, Rosser stated that since it was clear that $\text{Prov}(\varphi_m \leftrightarrow \exists x\ \{\varphi_x \text{ ProofOf } \varphi_n \wedge \neg\exists z\ (z \leq x \wedge \varphi_z \text{ ProofOf } \neg\varphi_n)\})$ is true then $\text{Prov}(\varphi_n) \rightarrow \text{Prov}(\varphi_m)$ is also true. However, $\varphi_m = \text{Prov}(\varphi_n) = \exists x\ (\varphi_x \text{ ProofOf } \varphi_n \wedge \neg\exists z\ \{z \leq x \wedge \varphi_z \text{ ProofOf } \neg\varphi_n\})$, due to the definitions of Rosser's provability predicate and proof relation, and thus the equivalence is obviously true by definition. As a result, since $\text{Prov}(\varphi_m) = \text{Prov}(\text{Prov}(\varphi_n))$ then it is also obvious that Rosser assumed that the truth of $\text{Prov}(\varphi_n)$ is sufficient to infer that a proof of $\text{Prov}(\varphi_n)$ exists.

The reasoning that Rosser employed to prove that $\text{Prov}(\neg\varphi_n) \rightarrow \text{Prov}(\neg\varphi_m)$ is true is analogous to the reasoning that he employed to prove that $\text{Prov}(\varphi_n) \rightarrow \text{Prov}(\varphi_m)$ is true. First, he proved that if a proof of $\neg\varphi_n$ exists then $\forall x\ (\neg\{\varphi_x \text{ ProofOf } \varphi_n\} \vee \exists z\ \{z \leq x \wedge \varphi_z \text{ ProofOf } \neg\varphi_n\})$ is true. Thus, since $\neg\text{Prov}(\varphi_n) = \forall x\ (\neg\{\varphi_x \text{ ProofOf } \varphi_n\} \vee \exists z\ \{z \leq x \wedge \varphi_z \text{ ProofOf } \neg\varphi_n\})$, due to the definitions of Rosser's provability predicate and proof relation as well as assuming $\neg\exists x\ \varphi(x) \rightarrow \forall x\ \neg\varphi(x)$ and $\neg\neg\varphi \rightarrow \varphi$, then $\text{Prov}(\neg\varphi_n) \rightarrow \neg\text{Prov}(\varphi_n)$ is true. Hence, he then stated that $\text{Prov}(\neg\varphi_n) \rightarrow \text{Prov}(\neg\varphi_m)$ is also true. As a result, since $\text{Prov}(\neg\varphi_m) = \text{Prov}(\neg\text{Prov}(\varphi_n))$ then it is again obvious that Rosser assumed that the truth of $\neg\text{Prov}(\varphi_n)$ is sufficient to infer that a proof of $\neg\text{Prov}(\varphi_n)$ exists.

Therefore, the principal assumption underlying the concept of validity that Rosser employed in his proof is, for any formula $\varphi(x)$, both $\varphi(x) \rightarrow \text{Prov}(\varphi(x))$ and $\neg\varphi(x) \rightarrow \text{Prov}(\neg\varphi(x))$ are true, which Rosser assumed without proof or any explicit justification. This assumption is also evident in the preceding proof of theorem 1 within the statement of the assumptions in 4 and 6. Although this assumption is similar to what Gödel proved with his theorem 5, the assumption cannot be justified by Gödel's theorem 5 because the relation in the antecedent of the theorem is required to be primitive recursive and neither $\text{Prov}(\varphi_n) = \exists x\ (\varphi_x \text{ ProofOf } \varphi_n \wedge \neg\exists z\ \{z \leq x \wedge \varphi_z \text{ ProofOf } \neg\varphi_n\})$ or $\neg\text{Prov}(\varphi_n) = \forall x\ (\neg\{\varphi_x \text{ ProofOf } \varphi_n\} \vee \exists z\ \{z \leq x \wedge \varphi_z \text{ ProofOf } \neg\varphi_n\})$ is primitive recursive due to the unbounded quantifiers within their definitions. Furthermore, Rosser rejected primitive recursion as a basis of his proof because he stated that his proof, which is based on recursive and recursively enumerable sets, is more general than Gödel's proof, which is based on primitive recursion and proof-theoretic deductive proofs.

As a consequence, the crucial difference between Gödel's proof and Rosser's proof is that Gödel's proof may be considered to be constructive whereas Rosser's proof contains assumptions that are not generally considered to be constructive. However, this can only be stated informally because a formal analogue of what is generally considered to be constructive does not exist and thus it cannot be formally proven.

Finally, it is of interest to explicitly note that the Rosser sentence that Rosser employed in his proof is primitive recursive and thus it is considered to be constructive whereas the definition of simple consistency employed in the preceding proof of theorem 1 is not primitive recursive and thus it is not considered to be constructive. However, since Rosser's proof contains additional assumptions that are not generally considered to be constructive then the fact that the definition of simple consistency is not primitive recursive is immaterial to the primary objective of the theorem.

## 10. TARSKI, MOSTOWSKI AND ROBINSON'S PROOF

Tarski, Mostowski and Robinson constructed a metamathematical proof that they claim is simpler and more general than Gödel's proof ([5] pg. 46 – 47). Instead of requiring the formal system in their proof to be capable of defining the entire class of primitive recursive functions, they simply require the formal system to be capable of constructing a bijective function between the formulas in the language of the system and the natural numbers. Hence, intuitively, they merely require an implementation of the metamathematical notion of the arithmetization of syntax. And instead of requiring the formal system in their proof to be capable of defining the proof-theoretic concept of deductive proof, which is the concept of provability employed in Gödel's proof, they merely require the informal and thus purportedly more general concept of validity within the formal system. Hence, intuitively, the metamathematical concept of validity employed in their proof accommodates any definition of what constitutes a valid formula within the formal system.

The basis of their proof is the notion of definability within the formal system, which we will now describe. For any function $f\colon \mathbb{N}^n \to \mathbb{N}$, where $\mathbb{N}$ denotes the set of natural numbers and $\mathbb{N}^n$ denotes the set of $n$-tuples of natural numbers, the function $f$ is definable within the formal system if some formula $\varphi(x_1, \ldots, x_n, y)$ exists within the formal system such that $f(x_1, \ldots, x_n) = y \leftrightarrow \varphi(x_1, \ldots, x_n, y)$. And then for any set $S$, the set $S$ is definable within the formal system if some formula $\varphi(x)$ exists within the formal system such that $s \in S \leftrightarrow \varphi(s)$ and thus also $s \notin S \leftrightarrow \neg\varphi(s)$.

Finally, $T$ will be employed to denote any formal system that meets the preceding requirements, which is any formal system within which the arithmetization of syntax and any notion of validity within the formal system can be formalized. In addition, as in Gödel's proof $[\varphi]$ will be employed to denote the unique natural number that is the Gödel number of the formula $\varphi$ in $T$ and $\varphi_n$ will be employed to denote the formula in $T$ whose Gödel number is $n$ and thus $[\varphi_n] = n$.

**Tarski, Mostowski and Robinson:** within any formal system $T$ that is simply consistent, the diagonal function $f_d(n) = [\varphi_n(n)]$ and the set $V$ that contains the Gödel numbers of all the valid sentences in $T$ are not both definable within $T$.

*Proof*: assume the contrary, which is that the diagonal function $f_d\colon \mathbb{N} \to \mathbb{N}$ and the set $V$ are both definable within the formal system $T$. As a result, the formulas $\Phi(x, y)$ and $\Psi(x)$, which define the diagonal function and the set $V$ respectively, must both exist within $T$. Therefore:

**1:** $f_d(x) = y \leftrightarrow \Phi(x, y)$

**2:** $n \in V \leftrightarrow \Psi(n)$

**3:** $n \notin V \leftrightarrow \neg\Psi(n)$

Next, the formula $\varphi_m(x)$ is defined:

**4:** $\varphi_m(x) \equiv_{\text{df}} \forall n\, (\Phi(x, n) \rightarrow \neg\Psi(n))$

Substituting $m = [\varphi_m(x)]$ for $x$ in 4:

**5:** $\varphi_m(m) = \forall n\, (\Phi(m, n) \rightarrow \neg\Psi(n))$

And then due to 1 with $x = m$:

**6:** $f_d(m) = y \leftrightarrow \Phi(m, y)$

And due to the definition of the diagonal function:

**7:** $f_d(m) = [\varphi_m(m)]$

If the sentence $\varphi_m(m)$ is true in $T$ then the sentence $\neg\Psi(f_d(m))$ is also true in $T$, due to 5, 6 and 7. Inversely, if the sentence $\varphi_m(m)$ is not true in $T$ then $\varphi_m(m)$ is not valid in $T$ and thus $f_d(m) \notin V$ is true in $T$, due to 7 and the definition of the set $V$, and then the sentence $\neg\Psi(f_d(m))$ is also true in $T$, due to 3. Therefore, in either case:

**8:** The sentence $\neg\Psi(f_d(m))$ is true within the formal system $T$.

However, the sentence $\forall n\, (\Phi(m, n) \rightarrow \neg\Psi(n))$ is true in $T$, due to 1, 7 and 8, and thus $\varphi_m(m)$ is true in $T$, due to 5. Since $\varphi_m(m)$ is true in $T$ then $\varphi_m(m)$ is valid in $T$ and thus $f_d(m) \in V$ is true in $T$, due to 7 and the definition of the set $V$, but then due to 2:

**9:** The sentence $\Psi(f_d(m))$ is true within the formal system $T$.

Consequently, due to 8 and 9, the formal system $T$ is simply inconsistent. ■

The proof that Tarski, Mostowski and Robinson constructed is inherently metamathematical because the proof does not contain explicit formal definitions of the arithmetization of syntax and the concept of validity within the formal system. Therefore, although it is claimed that the proof is simpler and more general than Gödel's proof, the proof crucially relies on the formal definitions constructed by Gödel in his proof to assert that the proof is not entirely theoretical. Hence, the proof is simpler only in the sense that it does not include those formal definitions. However, the proof is indeed more general than Gödel's proof in the sense that Gödel's proof is inherently proof-theoretic, due to Gödel's exclusive focus on deductive proofs, whereas this proof is applies to any formal system that is capable of defining any concept of validity within the formal system.

Although Tarski, Mostowski and Robinson's metamathematical proof cannot supplant Gödel's constructive proof, their proof may be considered to entirely supersede Rosser's proof because their proof is at least as general as Rosser's proof and may also be considered to be simpler. Gödel stated that his first incompleteness theorem applies to any formal system that satisfies three requirements ([1.a] pg. 181). These requirements are, (1) the primitive recursive functions can be defined within the formal system, (2) for any enumeration of the formulas in the language of the system, $n \in V \leftrightarrow \neg\mathrm{Prov}(\varphi_n(n))$ can be defined within the formal system, and finally (3) the formal system is ω-consistent. It has already been described how Tarski, Mostowski and Robinson's proof satisfies the first two requirements, but their proof, like Rosser's proof, apparently dispenses with the last requirement.

In conclusion, similar to the Gödel sentence $\varphi_p(p) = \forall n\ \neg(\varphi_n \text{ ProofOf } \varphi_p(p))$ and the Rosser sentence $\varphi_p(p) = \forall n\ \neg(\varphi_n \text{ ProofOf } \varphi_p(p) \wedge \neg\exists z\ \{z \leq n \wedge \varphi_z \text{ ProofOf } \neg\varphi_p(p)\})$ what may be called the Tarski, Mostowski and Robinson sentence $\varphi_m(m) = \forall n\ (\Phi(m, n) \rightarrow \neg\Psi(n))$ asserts that the sentence $\varphi_m(m)$ is not valid and thus the sentence states of itself that it is not true.

## 11. SUMMARY

In the preceding formal analysis we attempted to compose a comprehensive as well as intuitively meaningful explication of the proof of Gödel's first incompleteness theorem without violating one of the primary objectives of our proof by introducing interpretations or suppositions that cannot be directly justified by what is contained in Gödel's proof. However, some of the questions raised by Gödel's proof have never been formally resolved and thus they do not have a generally accepted solution. As a result, the preceding formal analysis is naturally and necessarily incomplete.

Therefore, along with the expected intuitive summary of the concept of incompleteness as well as some additional theorems that are immediately implied by Gödel's proof, the following summary will also include an examination of some of the questions raised by Gödel's proof that could not be directly formally addressed in the preceding formal analysis.

### 11.1 A Simple Intuitive Description of the Proofs of Incompleteness

We begin with a simple, semi-formal metamathematical proof that reveals the underlying structure of the proofs of incompleteness. Only an undefined provability predicate will be employed in the proof, which may be considered to represent any formal definition of validity within the system. Therefore, for any sentence φ in the language of some formal system:

**1:** $\varphi \rightarrow \mathrm{Prov}(\varphi)$

**2:** $\neg\varphi \rightarrow \mathrm{Prov}(\neg\varphi)$

1 and 2 simply state that if the sentence φ or ¬φ is true then φ or ¬φ is provable within the formal system respectively. Since φ denotes any sentence in the language of the system, we then define:

**3:** $\varphi \equiv_{\mathrm{df}} \neg\mathrm{Prov}(\varphi)$

Thus, φ is defined as the sentence that states of itself that it is not provable. Hence, due to 3 and then if $\neg\neg\varphi \rightarrow \varphi$ is assumed:

**4:** $\neg\varphi = \neg\neg\mathrm{Prov}(\varphi) = \mathrm{Prov}(\varphi)$

If φ is assumed to be true then Prov(φ) must be true, due to 1, but then ¬Prov(φ) must also be true, due to the definition of the sentence φ in 3, which is a contradiction if the system is assumed to be simply consistent. Subsequently, if ¬φ is assumed to be true then Prov(¬φ) must be true, due to 2, but then Prov(φ) must also be true, due to the equality in 4, which is a contradiction if the formal system is assumed to be simply consistent. The contradictions may be clearly seen by substituting the definition of φ in 3 for φ in the antecedent of 1 and 2 using the equality in 4.

**5:** ¬Prov(φ) → Prov(φ)

**6:** Prov(φ) → Prov(¬φ)

It is of interest to explicitly note that the contradiction inferred from 1 and 3, which is expressed in 5, requires the system, not the formal system, to be simply consistent. Although this reference may be interpreted as referring to the formal system, it may also be interpreted as referring to what could be called the metamathematical system, which is the theoretical system that expresses the naïve theory that is the foundation of the formal system.

As a result, the source of the contradictions as well as the underlying paradox becomes evident, however, to unduly emphasize the paradox underlying the proof of Gödel's first incompleteness theorem would violate Gödel's conceptualization of what his proof actually proved. This may be clearly seen in the following statement by Gödel ([8] pg. 179). In response to Wittgenstein's comments about Gödel's proof in *Remarks on the Foundation of Mathematics* [9], Gödel stated, "He [Wittgenstein] did *not* understand it (or pretended not to understand it). He interpreted it as a kind of logical paradox, while in fact it is just the opposite, namely a mathematical theorem within an absolutely uncontroversial part of mathematics (finitary number theory or combinatorics)."

In the preceding statement Gödel characterized the formal axiomatic system that he constructed in his proof as an absolutely uncontroversial part of mathematics because the axioms of his formal system consist of only well known and generally accepted axioms from elementary arithmetic and mathematical logic (Appendix B: Gödel's System). Gödel's definition of the primitive recursive functions, which defines the elementary operations of arithmetic within Gödel's formal system, is fundamentally based on the Peano Axioms (Axioms I). The list of relations that Gödel defined in his proof, which culminate in the definitions of Gödel's proof relation and provability predicate, are based on axioms from mathematical logic, specifically, axioms from Propositional Logic (Axioms II), Predicate Logic (Axioms III), and Set Theory (Axioms IV and V).

Gödel's characterization of his proof in the preceding statement then becomes clear when the principal objective of Gödel's proof is considered. The objective of Gödel's proof was stated in his second incompleteness theorem, which asserted that the formal system in Gödel's proof could not prove its own simple consistency and thus his proof was to some extent a response to Hilbert's program (section 2.2: Consistency, Completeness and Hilbert's Program). As a consequence, since a formal system is simply consistent if no sentence exists in the language of the system such that both the sentence and its negation can be proven to be true within the system then the necessity of stating and proving the Gödel sentence or some equivalent sentence becomes clear. Hence, according to Gödel, it is this intrinsic necessity that transforms the Gödel sentence from a statement that may be perceived to be a contrived absurdity that is concocted merely to construct a logical paradox into a meaningful statement.

The inadequacy of the preceding metamathematical proof and the logical paradox that it expresses may be clearly seen, at least in part, due to its inability to provide any insight into the role of the concept of ω-inconsistency within a proof of incompleteness.

### 11.2 Two Formal Interpretations of the Concept of ω-inconsistency

In [6.a] Tarski composed a formal set-theoretic definition of ω-consistency as well as what he called ω-completeness in an attempt to provide an intuitively meaningful characterization of the concept of ω-inconsistency. Therefore, we will employ a simple set-theoretic construction within a model-theoretic proof to illustrate the basic idea underlying Tarski's definition of ω-inconsistency.

First, we define the set of every finite subset of the natural numbers, which we denote by $\boldsymbol{P}_\omega(\mathbb{N})$, where $\mathbb{N}$ denotes the set of natural numbers. We then define the bijective function $f$: $\boldsymbol{P}_\omega(\mathbb{N}) \to \mathbb{N}$ as $f(s) = \sum_{n\in s} 2^n$, where $s \in \boldsymbol{P}_\omega(\mathbb{N})$, $n \in \mathbb{N}$ and $\sum_{n\in s}$ indicates the sum ranging over every $n \in s$. As is customary, we define the result of the empty summation to be 0 and thus $f(\varnothing) = 0$, where $\varnothing$ denotes the empty set or the set that does not contain any members. As a result, since $f$ is a bijective function then $f$ is a one-to-one correspondence between $\boldsymbol{P}_\omega(\mathbb{N})$ and $\mathbb{N}$ and thus $f$ is also a complete enumeration of every $s \in \boldsymbol{P}_\omega(\mathbb{N})$. Hence, if $s_n$ denotes the $s \in \boldsymbol{P}_\omega(\mathbb{N})$ such that $f(s) = n$ then $\boldsymbol{P}_\omega(\mathbb{N}) = \{s_0, s_1, s_2, \ldots\}$. Finally, we employ the customary notation $\cup$ to denote the standard set-theoretic operation of union and then also, for any $n \geq 0$, we will employ the notation $\cup^n s_i$ to denote the union of the sets $s_0, \ldots, s_n$ or $\cup^n s_i \equiv_{df} s_0 \cup \ldots \cup s_n$.

Therefore, we now employ a proof by mathematical induction to prove that the result of the union of every $s_0, s_1, s_2, \ldots \in \boldsymbol{P}_\omega(\mathbb{N})$ or $s_0 \cup s_1 \cup s_2 \cup \ldots$ is a finite set or a set that contains a finite number of members. The induction basis is vacuously true because $s_0 = \varnothing$ and $\varnothing$ is clearly a finite set. Next, to prove the induction step we assume the induction hypothesis or, for any $n \geq 0$, the result of $\cup^n s_i$ is a finite set and then we prove that the result of $\cup^{n+1} s_i$ must also be a finite set. Since $\cup^{n+1} s_i = (\cup^n s_i) \cup s_{n+1}$ then the result of $\cup^{n+1} s_i$ must be a finite set because the result of the union of any two finite sets is clearly a finite set and the result of $\cup^n s_i$ is a finite set by hypothesis and the set $s_{n+1}$ is a finite set by definition. As a consequence, it has been proven that the result of the union of every $s_0, s_1, s_2, \ldots \in \boldsymbol{P}_\omega(\mathbb{N})$ must be a finite set. However, this cannot be true because for every $n \in \mathbb{N}$ some $s \in \boldsymbol{P}_\omega(\mathbb{N})$ exists such that $n \in s$, due to the definition of $\boldsymbol{P}_\omega(\mathbb{N})$. Thus, the result of the union of every $s_0, s_1, s_2, \ldots \in \boldsymbol{P}_\omega(\mathbb{N})$ is $\mathbb{N}$ and $\mathbb{N}$ is obviously not a finite set.

A little reflection on the preceding proof by mathematical induction reveals the flaw in the proof. It was asserted that what was to be proven by mathematical induction was that the result of the union of every $s_0, s_1, s_2, \ldots \in \boldsymbol{P}_\omega(\mathbb{N})$ or $s_0 \cup s_1 \cup s_2 \cup \ldots$ is a finite set, but what was actually proven was that, for any $n \geq 0$, the result of the union of $s_0, \ldots, s_n \in \boldsymbol{P}_\omega(\mathbb{N})$ or $\cup^n s_i$ is a finite set. The difference between $s_0, \ldots, s_n \in \boldsymbol{P}_\omega(\mathbb{N})$ for any $n \geq 0$ and $s_0, s_1, s_2, \ldots \in \boldsymbol{P}_\omega(\mathbb{N})$ was characterized by Tarski as the difference between the concepts of a potential infinity and an actual or completed infinity ([6.a] pg. 293). Thus, more careful reflection on the preceding proof and its flaw reveals that if the concepts of a potential infinity and a completed infinity are not explicitly and clearly distinguished then the proof is not actually flawed, but rather the formal system within which the proof was constructed is ω-inconsistent, at least according to Tarski.

Within this context the concept of a potential infinity might be more accurately characterized as a constructive infinity that is related to the concept of evaluated truth and the concept of a completed infinity might be more accurately characterized as a definitional infinity that is related to the concept of asserted truth.

Subsequently, in [7.a] Quine defined an ω-inconsistent formal system as a system within which no comprehensive formal definition of the natural numbers exists in the language of the system. As a result, Quine called an ω-inconsistent formal system *numerically insegregative* ([7.a] pg. 118). A formal system that is numerically insegregative and thus also ω-inconsistent may be intuitively described as a formal system within which no formula exists in the language of the system that is true for every $n \in \mathbb{N}$ and only for every $n \in \mathbb{N}$.

Quine used the following example to describe a numerically insegregative formal system: assume that some formula $\varphi(x)$ exists in the language of the system such that $\varphi(x) \rightarrow x \in \mathbb{N}$. However, since the system is numerically insegregative then some $n$ exists such that $\neg\varphi(n) \wedge n \in \mathbb{N}$. Hence, $\varphi(x)$ cannot be purported to define the natural numbers within the system. Subsequently, assume that some other formula $\Psi(x)$ exists in the language of the system such that $\Psi(x) \rightarrow x \in \mathbb{N}$ and $\Psi(n)$ is true for the $n \in \mathbb{N}$ for which $\varphi(n)$ is false. As a result, the formula $\Phi(x) \equiv_{df} \varphi(x) \vee \Psi(x)$ exists in the language of the system and thus $\Phi(x) \rightarrow x \in \mathbb{N}$ and $\Phi(n)$ is true for the $n \in \mathbb{N}$ for which $\varphi(n)$ is false. However, since the system is numerically insegregative then some $m \neq n$ exists such that $\neg\Phi(m) \wedge m \in \mathbb{N}$. As a consequence, for any formula in the language of the formal system that may be purported to define the natural numbers within the system, some $n \in \mathbb{N}$ will always exist for which the formula is false since the system is numerically insegregative.

Hence, the primary focus of Tarski's definition is the inadequacy of any formal system to formally deal with the concepts of a potential infinity and a completed infinity whereas the primary focus of Quine's definition is the inadequacy of any specific formal system to formally define and prove theorems about the entire infinite set of natural numbers. However, both definitions are to some extent inadequate due to their inability to clearly and conclusively resolve some basic questions about the concept of ω-inconsistency within the context of Gödel's proof.

One such question is expressed by the inequality $\forall x\ \text{Prov}(\varphi(x)) \neq \text{Prov}(\forall x\ \varphi(x))$, which must be true due to Gödel's definition of ω-inconsistency. Gödel formally defined an ω-inconsistency as a formal system that may be simply consistent and within which $\forall x\ \text{Prov}(\varphi(x)) \wedge \text{Prov}(\neg\forall x\ \varphi(x))$ is true for some formula $\varphi(x)$ in the language of the system. Thus, if $\forall x\ \text{Prov}(\varphi(x)) = \text{Prov}(\forall x\ \varphi(x))$ were true then $\text{Prov}(\forall x\ \varphi(x)) \wedge \text{Prov}(\neg\forall x\ \varphi(x))$ would also be true in an ω-inconsistent system.

However, if the system is also simply consistent then Prov($\forall x\ \varphi(x)$) $\wedge$ Prov($\neg\forall x\ \varphi(x)$) is clearly a contradiction because the negation of the sentence $\forall x\ \varphi(x)$ is the sentence $\neg\forall x\ \varphi(x)$. As a result, if $\forall x$ Prov($\varphi(x)$) = Prov($\forall x\ \varphi(x)$) were true then an ω-inconsistent formal system could never also be simply consistent, which contradicts Gödel's definition of ω-inconsistency.

Tarski's definition of ω-inconsistency suggests that $\forall x$ Prov($\varphi(x)$) signifies a potential infinity and Prov($\forall x\ \varphi(x)$) signifies an actual or completed infinity. However, the concepts of a potential infinity and a completed infinity remain mostly intuitive and thus also entirely theoretical. Thus, Tarski's definition merely gives $\forall x$ Prov($\varphi(x)$) $\neq$ Prov($\forall x\ \varphi(x)$) an intuitive interpretation and is thus is incapable of formally analyzing or resolving the inequality within the context of Gödel's proof. Tarski's *rule of infinite induction* ([6.a] pg. 294), which is his proposed resolution of the inequality, reveals the inadequacy of Tarski's definition within the context of Gödel's proof. In essence, Tarski's rule of infinite induction simply asserts that $\forall x$ Prov($\varphi(x)$) $\rightarrow$ Prov($\forall x\ \varphi(x)$), but this has already been shown to contradict Gödel's definition of ω-inconsistency.

Quine's definition of ω-inconsistency suggests that the language of an ω-inconsistent formal system is effectively incapable of defining a formula $\varphi(x)$ that is true if and only if the value of the variable $x$ is a natural number. As a result, an ω-inconsistent system is effectively incapable of defining the infinite set $\mathbb{N}$ and is thus numerically insegregative. Therefore, Prov($\forall x\ \varphi(x)$) can be true in an ω-inconsistent system, but $\forall n \in \mathbb{N}$ Prov($\varphi(n)$) can never be true because $\mathbb{N}$ cannot be defined within the system, which explains how $\forall x$ Prov($\varphi(x)$) $\neq$ Prov($\forall x\ \varphi(x)$) can be true when the values of the variable $x$ range over the natural numbers. However, $\mathbb{N}$ can be defined within the formal system in Gödel's proof because axioms I (Appendix B: Gödel's System) are based on the Peano Axioms, which are the *de facto* standard definition of $\mathbb{N}$. Thus, Quine's definition is also inadequate within the context of Gödel's proof.

### 11.3 An Intuitionistic Interpretation of Gödel's Concept of ω-inconsistency

The intuitive foundation of Gödel's concept of ω-inconsistency can be related to the fundamental difference between the constructive and definitional interpretations of the existential quantifier (section 6: The Intuitionistic Acceptability of Gödel's Proof). The constructive interpretation of the existential quantifier essentially states that, for any formula $\varphi(x)$, the sentence $\exists x\ \varphi(x)$ asserts that some specific defined value of the variable $x$ exists such that $\varphi(x)$ can be proven to be true for this value, in symbols: $\exists x\ \varphi(x) \rightarrow \varphi(c)$, where $c$ is some fixed value. And then the definitional interpretation of the existential quantifier simply states that some value of the variable $x$ exists for which $\varphi(x)$ is true, although no such actual value may be known.

If a proof of the sentence $\exists x\ \varphi(x)$ exists for some formula $\varphi(x)$ in the language of a formal system, but a specific defined value of the variable $x$ is not explicitly identified in the proof for which $\varphi(x)$ can be proven to be true, then the definitional interpretation of the existential quantifier has been assumed. If a specific defined value for which $\varphi(x)$ is true is never identified then this suggests that no such value actually exists and thus $\exists x\ \varphi(x)$ is not actually true. However, a specific defined value for which $\varphi(x)$ is true may actually exist although no such value is ever identified for whatever reason. Regardless, if the formal system within which the proof of the sentence $\exists x\ \varphi(x)$ exists is simply consistent then the mere existence of the proof asserts that a proof of the sentence $\neg\exists x\ \varphi(x)$ cannot exist and neither can a proof of the sentence $\forall x\ \neg\varphi(x)$ if $\forall x\ \neg\varphi(x) \rightarrow \neg\exists x\ \varphi(x)$ is true. As a result, since a proof of the sentence $\forall x\ \neg\varphi(x)$ cannot exist then the sentence $\neg\forall x\ \neg\varphi(x)$ can be inferred if $\varphi \vee \neg\varphi$ (the LEM) is true. Conversely, if $\exists x\ \varphi(x) \rightarrow \neg\forall x\ \neg\varphi(x)$ is true, which is justified by Gödel's relation 32.d (Appendix B: Gödel's System), then the sentence $\neg\forall x\ \neg\varphi(x)$ can be directly inferred from the proof of the sentence $\exists x\ \varphi(x)$. In either case, if a specific defined value of the variable $x$ has never been identified for which $\varphi(x)$ can be proven to be true then the sentence $\exists x\ \varphi(x)$ does not necessarily assert that a specific defined value exists for which $\varphi(x)$ can be proven to be true and neither does the sentence $\neg\forall x\ \neg\varphi(x)$, in symbols: $\neg\forall x\ \neg\varphi(x) \rightarrow \varphi(c)$ is not necessarily true. As a result, even though $\neg\forall x\ \neg\varphi(x)$ is taken to be true $\varphi(c)$ is not necessarily true for any fixed value $c$.

Therefore, within Gödel's proof (section 4: Gödel's Proof), Prov($\neg\forall x\ \varphi_r(x)$) may be interpreted as stating that the existence of a deductive proof of the sentence $\neg\forall x\ \varphi_r(x)$ within the formal system does not then necessarily assert that $\neg\forall x\ \varphi_r(x) \rightarrow \neg\varphi_r(c)$ is true for some fixed value $c$. Indeed, no specific defined value of the variable $x$ was identified in Gödel's proof for which $\neg\varphi_r(x)$ can be proven to be true. As a consequence, since $\forall x$ Prov($\varphi_r(x)$) is interpreted as stating that, for every fixed value $c$, a deductive proof of the sentence $\varphi_r(c)$ exists within the formal system then the sentence $\text{Prov}(\neg\forall x\ \varphi_r(x)) \wedge \forall x\ \text{Prov}(\varphi_r(x))$, which expresses the concept of ω-inconsistency within Gödel's proof, does not necessarily entail a simple inconsistency.

Subsequently, we assert that unlike either Tarski's or Quine's definition of ω-inconsistency this interpretation of ω-inconsistency is capable of clearly expressing the precise difference between $\forall x\ \text{Prov}(\varphi(x))$ and $\text{Prov}(\forall x\ \varphi(x))$ and thus also why the inequality $\forall x\ \text{Prov}(\varphi(x)) \neq \text{Prov}(\forall x\ \varphi(x))$ may be true (section 11.2: Two Formal Interpretations of the Concept of ω-inconsistency).

Intuitively, for any formula $\varphi(x)$, $\text{Prov}(\forall x\ \varphi(x))$ states that a single deductive proof exists within the formal system of the sentence that states that $\varphi(x)$ is true for every value of the variable $x$ and $\forall x\ \text{Prov}(\varphi(x))$ states that, for every value of the variable $x$, a deductive proof of $\varphi(x)$ exists within the formal system for that individual value of the variable $x$. Therefore, it appears to be clear that both $\text{Prov}(\forall x\ \varphi(x))$ and $\forall x\ \text{Prov}(\varphi(x))$ reflect the same intuitive condition, which is that $\varphi(x)$ is true within the formal system for every possible value of the variable $x$.

Formally, for any numeric formula $\varphi(x)$, it is clear that it can be inferred from the axioms of Gödel's formal system that $\text{Prov}(\forall x\ \varphi(x))$ is the result of proof that depends on the principle of mathematical induction because the only axiom in the system that will permit a deductive proof of the sentence $\forall x\ \varphi(x)$ to exist within the system is axiom I.3 (Appendix B: Gödel's System), which is a set-theoretic formalization of the principle of mathematical induction. And, for any primitive recursive numeric formula $\varphi(x)$, it can be inferred from the axioms of Gödel's formal system along with Gödel's theorem 5 that $\forall x\ \text{Prov}(\varphi(x))$ is the result of a proof by mathematical induction (section 5: The Concepts of ω-inconsistency and Mathematical Induction). Therefore, it appears to be clear that both $\text{Prov}(\forall x\ \varphi(x))$ and $\forall x\ \text{Prov}(\varphi(x))$ are, in essence, the same formal proof because they are both formal analogues of the metamathematical principle of mathematical induction.

Therefore, in either case and within the context of Gödel's proof, for any numeric formula $\varphi(x)$, where the values of variable $x$ range over the natural numbers, $\text{Prov}(\forall x\ \varphi(x))$ represents a formal deductive proof of the sentence $\forall x\ \varphi(x)$ from the premises: axiom I.3, $\varphi(0)$, and $\varphi(n) \rightarrow \varphi(S(n))$, where $n$ denotes an arbitrary natural number and $S$ denotes the successor operation. And, for any primitive recursive numeric formula $\varphi(x)$, $\forall x\ \text{Prov}(\varphi(x))$ represents a metamathematical proof of $\forall x\ \varphi(x)$ that is justified by the premise that a formal deductive proof of $\varphi(x)$ must exist whenever the value of the variable $x$ is any natural number. The proof is called metamathematical because its premise is a clear expression of the customary metamathematical interpretation of the result of a proof by mathematical induction, which states that the numeric predicate or formula that is the focus of the proof must be true for every possible natural number. Furthermore, the proof is clearly metamathematical because $\forall x\ \text{Prov}(\varphi_r(x))$ does not represent a single formal deductive proof, but rather represents an infinite number of formal deductive proofs.

As a result, although it appears to be clear that both $\text{Prov}(\forall x\ \varphi(x))$ and $\forall x\ \text{Prov}(\varphi(x))$ reflect the same intuitive condition as well as the same formal proof, $\text{Prov}(\forall x\ \varphi(x))$ indicates that a formal deductive proof of the sentence $\forall x\ \varphi(x)$ exists within the formal system, whereas $\forall x\ \text{Prov}(\varphi(x))$ indicates that a metamathematical proof of $\forall x\ \varphi(x)$ exists that states that a formal deductive proof of $\varphi(x)$ exists within the formal system for every individual value of the variable $x$. Therefore, the inequality $\forall x\ \text{Prov}(\varphi(x)) \neq \text{Prov}(\forall x\ \varphi(x))$ is true because each side of the inequality indicates a proof of $\forall x\ \varphi(x)$ that is fundamentally and meaningfully different or is at least as fundamentally and meaningfully different as the constructive and definitional interpretations of the existential quantifier are considered to be.

### 11.4 The Concept of ω-inconsistency

The ω in ω-inconsistency refers to the transfinite ordinal number of the infinite set $\mathbb{N}$, which is the set that contains every natural number. The conceptualization underlying the ordinal ω is reflected in its set-theoretic construction, which characterizes the ordinal ω and thus also $\mathbb{N}$ as completed infinite totalities that consist of an infinite number of discrete members. This metamathematical conceptualization is in contrast to the constructive conceptualization of $\mathbb{N}$, which characterizes the natural numbers as being constructed by a perpetually unfolding process that never terminates and is thus never complete. Hence, the concept of ω-inconsistency that Gödel introduced and defined in his proof may be interpreted as challenging the customary set-theoretic assumption that an infinite set is nothing more than the sum or aggregate of its infinite discrete parts, which is the assumption that is intimated by the ordinal ω and also by the cumulative hierarchy in Set Theory.

The term inconsistency in ω-inconsistency refers to the concept of logical inconsistency, which is generally represented by its formal analogue of simple inconsistency. A formal system is defined to be simply consistent if no sentence exists in the language of the system such that both the sentence and its negation can be proven to be true within the formal system whereas a system is logically consistent if no two mutually absurd or contradictory statements exist within the system and both of the statements are true within the system. Hence, Gödel's concept of ω-inconsistency may be interpreted as postulating the possibility of the existence of a proof within the system that states that some formula or property is true of every individual natural number along with the existence of another proof within the system that states that the same formula or property cannot be true of every possible natural number.

Therefore, the conceptual difference between the two statements or conditions that constitute an ω-inconsistency in Gödel's proof may be stated clearly and precisely. The statement that a formula or property is true of every individual natural number is a positive statement about what is true of each specific defined natural number whereas the statement that a formula or property is not true of every possible natural number is a nonspecific negative statement about what cannot be true of every natural number in general. Hence, Gödel's concept of ω-inconsistency may be interpreted as postulating the possibility that, for some formula $\varphi(x)$ in the language of a formal system, a proof exists that states that $\varphi(x)$ must be true whenever the value of the variable $x$ is any specific defined natural number and a proof also exists that states that the sentence $\forall x\ \varphi(x)$ cannot be true when the values of the variable $x$ range over every possible natural number.

As a result, a practical interpretation of an ω-inconsistent formal system may be postulated to be, for some formula $\varphi(x)$ in the language of the system, a proof exists that employs the principle of mathematical induction and states that, for any specific defined $n \in \mathbb{N}$, the sentence $\varphi(n)$ must be true and a proof also exists that employs the principle of *reductio ad absurdum* and states that the sentence $\forall x\ \varphi(x)$ cannot be true when the values of the variable $x$ range over all of the natural numbers. However, this interpretation is not justified by what is contained in Gödel's proof because the proof-theoretic formal system in Gödel's proof does not contain a formalization of a method of proof that would permit a proof that employs the principle of *reductio ad absurdum* to exist within the formal system. Although this is in contrast to the fact that the proof of Gödel's first incompleteness theorem is itself a proof by contradiction, which is the customary formal analogue of a proof that employs the principle of *reductio ad absurdum*.

Inversely, a meaningful interpretation of an ω-inconsistent formal system may be postulated to be, for some formula $\varphi(x)$ in the language of the system, a proof can be composed that states that $\varphi(x)$ must be true within the formal system whenever any possible value is substituted for its free variable and another proof can be composed that states that the formal sentence $\forall x\ \varphi(x)$ cannot be true within the formal system or equivalently, if the LEM is assumed to be true, that the formal sentence $\neg\forall x\ \varphi(x)$ is true within the formal system. Although this interpretation can be justified by what is contained in Gödel's proof the two proofs can easily be perceived to contradict each other, at least intuitively, and thus the interpretation becomes dubious to some extent. Hence, the interpretation cannot actually be purported to be meaningful.

To illustrate the basic intuitive conceptualization underlying Gödel's concept of ω-inconsistency, which is devoid of any interpretation, we will employ a set-theoretic definition of ω-inconsistency that is essentially equivalent to the definition in Gödel's proof. First, we define the theory of a formal system as, $Th = \{n \in \mathbb{N} \mid \text{Prov}(\varphi_n)\}$. However, in this definition the provability predicate represents any means by which a sentence in the formal system can be determined to be true. As a result, the set $Th$ contains every natural number that is the Gödel number of a sentence in the language of the formal system that is true. Subsequently, for some formula $\varphi(x)$ in the language of a formal system, we define $k = [\forall x\ \varphi(x)]$ and also $k_n = [\varphi(n)]$ for every $n \in \mathbb{N}$, where $k \in \mathbb{N}$ and each $k_n \in \mathbb{N}$. An ω-inconsistency is then defined as $\forall n \in \mathbb{N}\ (k_n \in Th) \wedge [\neg\varphi_k] \in Th$, where $\neg\varphi_k$ denotes the sentence that is the negation of the sentence $\varphi_k = \forall x\ \varphi(x)$ and thus $[\neg\varphi_k] \in \mathbb{N}$ is the Gödel number of the sentence $\neg\varphi_k = \neg\forall x\ \varphi(x)$.

As a consequence, although $k \in Th \rightarrow \forall n \in \mathbb{N}\ (k_n \in Th)$ is clearly true due to Gödel's axiom III.1 (Appendix B: Gödel's System), Gödel's concept of ω-inconsistency postulates the possibility that $\forall n \in \mathbb{N}\ (k_n \in Th) \rightarrow k \in Th$ is not necessarily true within a simply consistent formal system. If $\forall n \in \mathbb{N}\ (k_n \in Th) \rightarrow k = [\varphi_k] \in Th$ were necessarily true then $\forall n \in \mathbb{N}\ (k_n \in Th) \wedge [\neg\varphi_k] \in Th$ could never be true within a simply consistent formal system and thus an ω-inconsistent formal system could never also be simply consistent. However, an ω-inconsistent formal system is not necessarily simply inconsistent according to the definition of ω-inconsistency in Gödel's proof. Therefore, the basic intuitive conceptualization underlying Gödel's concept of ω-inconsistency is that $\forall n \in \mathbb{N}\ (k_n \in Th) \rightarrow k \in Th$ and thus also $\forall n \in \mathbb{N}\ (k_n \in Th) \leftrightarrow k \in Th$ are not necessarily true within a simply consistent formal system.

In conclusion, the following statement by Gödel alludes to the paramount importance of the subtle distinction that we have attempted to exhibit ([8] pg. 78), "I believe that the true meaning of the opposition between things and concepts or between factual and conceptual truth is not yet completely understood in contemporary philosophy, but so much at least is clear: that in both cases one is faced with "solid facts", which are entirely outside the reach of arbitrary decisions."

### 11.5 Simple Intuitive Descriptions of ω-inconsistency and the Proofs of Incompleteness

We conclude our examination of the concept of ω-inconsistency with another simple, semi-formal metamathematical proof that reveals the underlying structure of Gödel's proof of incompleteness as well as the role of ω-inconsistency in his proof. This proof is a simple extension of our previous proof (section 11.1: A Simple Intuitive Description of the Proofs of Incompleteness).

An undefined proof relation will be employed in the proof, which like our previous proof may be considered to represent any formal definition of validity within the system. And then a provability predicate will also be employed that is defined in terms of the proof relation in the usual manner. In addition, some ability to reason over every sentence as well as every proof in the language of the formal system is necessary and thus an implementation of the metamathematical notion of the arithmetization of syntax is also required. Unlike in our previous proof, which only required some specific sentence to be defined, a formula with a single free variable replaces that sentence in this proof. Therefore, for any formula $\varphi(x)$ in the language of some formal system:

**1:** $\varphi(x) \rightarrow \text{Prov}(\varphi(x))$

**2:** $\neg\varphi(x) \rightarrow \text{Prov}(\neg\varphi(x))$

1 and 2 state that for any fixed value $c$ that is a valid substitution for the free variable in $\varphi(x)$ if the sentence $\varphi(c)$ or $\neg\varphi(c)$ is true then $\varphi(c)$ or $\neg\varphi(c)$ is provable within the formal system respectively. Since $\varphi(x)$ denotes any formula in the language of the system, we then define:

**3:** $\varphi(x) \equiv_{\text{df}} \neg(\varphi_x \text{ ProofOf } \forall n\ \varphi(n))$

Unlike in the previous proof, which defined a sentence that states of itself that it is not provable, in this proof the formula $\varphi(x)$ states that whenever any fixed value is substituted for its free variable then that value is not the Gödel number of a proof of the universal closure of itself, which is the sentence $\forall n\ \varphi(n)$. Hence, due to 3 and then if $\neg\neg\varphi \rightarrow \varphi$ is assumed:

**4:** $\neg\varphi(x) = \neg\neg(\varphi_x \text{ ProofOf } \forall n\ \varphi(n)) = \varphi_x \text{ ProofOf } \forall n\ \varphi(n)$

If it is assumed that a proof of the sentence $\forall n\ \varphi(n)$ does not exist then a fixed value $c = [\varphi_c]$ also does not exist such that $\varphi_c$ is a proof of the sentence $\forall n\ \varphi(n)$. As a result, a proof of the sentence $\varphi(c)$ must exist for every fixed value $c$ that is a valid substitution for the free variable in $\varphi(x)$, due to 1 and the definition of $\varphi(x)$ in 3, but this contradicts the assumption that a proof of the sentence $\forall n\ \varphi(n)$ does not exist if the formal system is assumed to be ω-consistent. Subsequently, if it is assumed that a proof of the sentence $\forall n\ \varphi(n)$ does exist then a fixed value $c = [\varphi_c]$ also exists such that $\varphi_c$ is a proof of the sentence $\forall n\ \varphi(n)$. As a result, a proof of the sentence $\neg\varphi(c)$ must exist for some fixed value $c$ that is a valid substitution for the free variable in $\varphi(x)$, due to 2 and the equality in 4, but this contradicts the assumption that a proof of the sentence $\forall n\ \varphi(n)$ exists if the formal system is assumed to be simply consistent. The contradictions may be more clearly expressed by substituting the definition of $\varphi(x)$ in 3 for $\varphi(x)$ in the antecedent of 1 and 2 using the equality in 4.

**5:** $\neg(\varphi_x \text{ ProofOf } \forall n\ \varphi(n)) \rightarrow \text{Prov}(\varphi(x))$

**6:** $\varphi_x \text{ ProofOf } \forall n\ \varphi(n) \rightarrow \text{Prov}(\neg\varphi(x))$

And then the contradictions may be clearly seen by adding the applicable quantifiers:

**7:** $\forall x\ \neg(\varphi_x \text{ ProofOf } \forall n\ \varphi(n)) \rightarrow \forall x\ \text{Prov}(\varphi(x))$

**8:** $\exists x\ (\varphi_x \text{ ProofOf } \forall n\ \varphi(n)) \rightarrow \exists x\ \text{Prov}(\neg\varphi(x))$

As a result, the source of the contradictions in Gödel's proof is clearly illustrated and thus the role of ω-inconsistency within his proof becomes evident as well.

Since this proof reflects the underlying structure of a proof of incompleteness that includes the concept of ω-inconsistency then our previous proof reflects the underlying structure of a proof of incompleteness that does not include the concept of ω-inconsistency. Hence, this proof reveals the underlying structure of Gödel's proof of incompleteness whereas our previous proof reveals the underlying structure of all of the other proofs of incompleteness that do not require the assumption that the formal system is ω-consistent, which includes Rosser's proof of incompleteness as well as Tarski, Mostowski and Robinson's proof of incompleteness.

Therefore, we will employ our previous proof to reveal the underlying structure of Rosser's proof of incompleteness. The sentence $\varphi$ in the previous proof will be defined to be the Gödel sentence $\varphi_p(p) = \forall n\ \varphi_r(n)$ in Gödel's proof (section 4: Gödel's Proof), in symbols: $\varphi \equiv_{\text{df}} \forall n\ \varphi_r(n)$, and then $\forall n\ \varphi_r(n)$ is substituted for $\varphi$ in 1 and 2 in our previous proof:

**1:** $\forall n\ \varphi_r(n) \rightarrow \text{Prov}(\forall n\ \varphi_r(n))$

**2:** $\neg\forall n\ \varphi_r(n) \rightarrow \text{Prov}(\neg\forall n\ \varphi_r(n))$

Due to the definition of the formula $\varphi_r$, then the definition of the relation $Q$, and then the equality $\varphi_p(p) = \forall n\ \varphi_r(n)$ along with a change of the bound variable for clarity:

**3:** $\forall n\ \varphi_r(n) = \forall n\ \varphi_q(n, p) = \forall n\ \neg(\varphi_n \text{ ProofOf } \varphi_p(p)) = \forall x\ \neg(\varphi_x \text{ ProofOf } \forall n\ \varphi_r(n))$

Hence, due to 3 and then Gödel's relation 32.d (Appendix B: Gödel's System):

**4:** $\neg\forall n\ \varphi_r(n) = \neg\forall x\ \neg(\varphi_x \text{ ProofOf } \forall n\ \varphi_r(n)) = \exists x\ (\varphi_x \text{ ProofOf } \forall n\ \varphi_r(n))$

If $\forall n\ \varphi_r(n)$ is assumed to be true then $\text{Prov}(\forall n\ \varphi_r(n))$ must also be true, due to 1. However, since $\text{Prov}(\forall n\ \varphi_r(n)) = \exists x\ (\varphi_x \text{ ProofOf } \forall n\ \varphi_r(n)) = \neg\forall x\ \neg(\varphi_x \text{ ProofOf } \forall n\ \varphi_r(n))$, due to the definition of the provability predicate and then Gödel's relation 32.d, if the system is assumed to be simply consistent then this contradicts the definition of the sentence $\forall n\ \varphi_r(n)$, due to the equality in 3. Subsequently, if $\neg\forall n\ \varphi_r(n)$ is assumed to be true then $\text{Prov}(\neg\forall n\ \varphi_r(n))$ must also be true, due to 2. However, since $\text{Prov}(\neg\forall n\ \varphi_r(n)) = \exists x\ (\varphi_x \text{ ProofOf } \neg\forall n\ \varphi_r(n))$, due to the definition of the provability predicate, if the formal system is assumed to be simply consistent then this contradicts the definition of the sentence $\neg\forall n\ \varphi_r(n)$, due to the equality in 4. The contradictions may be more clearly expressed by substituting the definition of $\forall n\ \varphi_r(n)$ for the sentence $\forall n\ \varphi_r(n)$ in the antecedent of 1 and 2 using the equalities in 3 and 4.

**5:** $\forall x\ \neg(\varphi_x \text{ ProofOf } \forall n\ \varphi_r(n)) \rightarrow \text{Prov}(\forall n\ \varphi_r(n))$

**6:** $\exists x\ (\varphi_x \text{ ProofOf } \forall n\ \varphi_r(n)) \rightarrow \text{Prov}(\neg\forall n\ \varphi_r(n))$

And the contradictions may then be clearly seen by employing the definition of the provability predicate as well as Gödel's relation 32.d:

**7:** $\forall x\ \neg(\varphi_x \text{ ProofOf } \forall n\ \varphi_r(n)) \rightarrow \neg\forall x\ \neg(\varphi_x \text{ ProofOf } \forall n\ \varphi_r(n))$

**8:** $\text{Prov}(\forall n\ \varphi_r(n)) \rightarrow \text{Prov}(\neg\forall n\ \varphi_r(n))$

As a result, it becomes evident that the basis of the contradictions in Rosser's proof is the logical paradox that is implicit within the self-referential and self-contradictory formal sentence that states of itself that it is not provable.

Our previous proof can also be employed in precisely the same manner to reveal the underlying structure of Tarski, Mostowski and Robinson's proof of incompleteness. First, the sentence $\varphi$ in our previous proof is defined to be the Tarski, Mostowski and Robinson sentence $\varphi_m(m)$, in symbols: $\varphi \equiv_{\text{df}} \varphi_m(m)$. Subsequently, 3 becomes $\varphi_m(m) = \forall n\ (n = [\varphi_m(m)] \rightarrow \neg\text{Prov}(\varphi_n))$ and then 4 becomes $\neg\varphi_m(m) = \neg\forall n\ (n = [\varphi_m(m)] \rightarrow \neg\text{Prov}(\varphi_n)) = \exists n\ (n = [\varphi_m(m)] \wedge \text{Prov}(\varphi_n))$, assuming $\neg\forall x\ \varphi(x) \rightarrow \exists x\ \neg\varphi(x)$ and $(\varphi \rightarrow \neg\Psi) \rightarrow \neg(\varphi \wedge \Psi)$ and thus also $\neg\neg\varphi \rightarrow \varphi$ are true. Since $n = [\varphi_n]$ then $\varphi_n$ denotes the sentence $\varphi_m(m)$, which is the sentence that results from substituting the value $m = [\varphi_m(x)]$ for the free variable in the formula $\varphi_m(x)$.

### 11.6 Some Additional Theorems About Incompleteness

We will now prove the following two theorems, which are immediately implied by the proof of Gödel's first incompleteness theorem. The customary definition of the Law of the Excluded Middle (LEM) as well as its intuitive interpretation will be assumed, which states that for any sentence in the language of a formal system either the sentence or its negation must be true, in symbols: $\varphi \vee \neg\varphi$. In addition, it will also be assumed that what is true within a formal system is synonymous with what can be proven within that formal system, which is essentially the intuitive corollary to proposition 1 (section 2.5: The Fundamental Concepts of Proof and Truth).

**Theorem 2:** the LEM does not necessarily hold in a formal system that is incomplete.

*Proof*: immediate due to the definition of an incomplete formal system in Gödel's proof, which states that a formal system is incomplete if, for any sentence in the language of the system, a proof does not exist within the system for neither the sentence nor for its negation ([1.a] pg. 147). ■

Thus, any intuitive or philosophical argument based on a semantic evaluation of truth that attempts to assert the validity or truth of the LEM is rejected as irrelevant because, in practice, the LEM clearly does not hold in a formal system that is incomplete. And then furthermore:

**Theorem 3:** a proof that employs the principle of *reductio ad absurdum* is not necessarily a valid method of proof in a formal system within which the LEM does not hold.

*Proof*: a proof that employs the principle of *reductio ad absurdum* will often assert that if either the sentence $\varphi$ or the sentence $\neg\varphi$ is proven to be absurd and thus not true then the other sentence must necessarily be true. Therefore, if a proof does not exist within a formal system for neither $\varphi$ or $\neg\varphi$, which must be the case for some sentences $\varphi$ and $\neg\varphi$ in the language of a system within which the LEM does not hold, then a proof that either $\varphi$ or $\neg\varphi$ is absurd and thus not true cannot subsequently assert that the other sentence must necessarily be true because both sentences may not be true within the formal system. ■

**Corollary 1:** $\neg\neg\varphi \rightarrow \varphi$ does not hold in a formal system within which the LEM does not hold.

*Proof*: immediate due to the proof of theorem 3, assuming $\neg\neg\varphi \equiv_{df} \neg\varphi$ *is False*. The corollary may also be inferred from Kleene's remark that $\neg\neg\varphi \rightarrow \varphi$ and $\varphi \vee \neg\varphi$ are interchangeable as axioms in the classical system ([4] pg. 120). ■

It should be emphasized that theorem 2 does not state that the LEM is false, but rather that it is not necessarily true. Theorem 3 then makes an equivalent statement about the validity of a proof that employs the principle of *reductio ad absurdum*. Hence, the preceding theorems may be interpreted as a justification of the intuitionistic assertion that the LEM is not necessarily valid in every case.

### REFERENCES

[1] Kurt Gödel Collected Works, vol. I publications 1929 – 1936 eds. Feferman, Solomon et al., Oxford: Oxford University Press 1986 ISBN: 978-0-19-514720-9
[a] Gödel, Kurt On formally undecidable propositions of *Principia mathematica* and related systems I, pg. 145 – 195
[2] From Frege to Gödel, A Source Book in Mathematical Logic, 1897 – 193 ed. van Heijenoort, Jean, Cambridge Mass.: Harvard University Press 1967 ISBN: 0-674-32449-8
[a] Gödel, Kurt On formally undecidable propositions of *Principia mathematica* and related systems I, pg. 596 – 616
[b] Hilbert, David On the infinite, pg. 369 - 392
[3] The Undecidable: Basic Papers on Undecidable Propositions, Unsolvable Problems and Computable Functions ed. Davis, Martin, New York: Dover Publications reprint 2004 ISBN: 0-486-43228-9
[a] Gödel, Kurt On formally undecidable propositions of *Principia mathematica* and related systems I, pg. 5 – 38
[b] Rosser, J. B. Extensions of some theorems of Gödel and Church, pg. 231 – 235
[4] Kleene, Stephen C. Introduction to Metamathematics, New York: Ishi Press International 2009 ISBN: 0-923891-57-9
[5] Tarski, Alfred; Mostowski, Andrzej; Robinson, Raphael M. Undecidable Theories, Studies in Logic and the Foundation of Mathematics, New York: Dover Publications reprint 2010 ISBN: 978-0-486-47703-9
[6] Logic, Semantics, Metamathematics ed. Corcoran, John, Indianapolis In.: Hackett Publishing Company 1983 ISBN: 978-0-915144-75-4
[a] Tarski, Alfred. Some observations on the concepts of ω-consistency and ω-completeness, pg. 279 – 295
[7] Selected Logic Papers Enlarged Edition ed. Quine, Willard Van Orman, Cambridge Mass.: Harvard University Press 1996 ISBN: 0-674-79837-6
[a] Quine, W. V. On ω-inconsistency and a so-called axiom of infinity, pg. 114 – 120
[8] Wang, Hao A Logical Journey, From Gödel to Philosophy, Cambridge Mass.: The MIT Press 1996 ISBN: 978-0-262-23189-3
[9] Wittgenstein, Ludwig Remarks on the Foundations of Mathematics, Cambridge Mass.: The MIT Press 1983 ISBN: 978-0-262-73067-9

## APPENDIX A: GÖDEL NUMBERS AND THE ARITHMETIZATION OF SYNTAX

We will formally define a process of constructing Gödel numbers that is equivalent to the process outlined informally in section 2.4: Gödel Numbers and the Arithmetization of Syntax. In the following definitions $\mathbb{N}$ denotes the set of natural numbers, $L$ denotes the language of some formal system, which consists of every finite sequence of symbols constructed from an alphabet of $n \geq 1$ unique defined symbols, and $L^*$ denotes the set of every finite ordered subset of $L$. In addition, a consecutive natural number is assigned to each of the symbols in the alphabet of $L$ starting with 1.

**A Method of Constructing Gödel Numbers:** the bijective functions $f_n$: $L \rightarrow \mathbb{N}$ and $g_n$: $L^* \rightarrow \mathbb{N}$ are defined as:

1. $f_n(s) = \sum_{c_i \in s} n^{i-1}(c_i)$
2. $g_n(S) = \prod_{s_i \in S} p_i^{f_n(s_i)}$

As is customary, for any $n > 0$ we define $n^0 = 1$ and then we also define the result of an empty summation as well as an empty multiplication to be 0.

In definition 1, $s \in L$ denotes some string, which is some ordered sequence of characters and a character is a symbol in the alphabet of $L$ that is not necessarily unique in $s$. Hence, $s = c_1, \ldots, c_m$, where $m \geq 0$ is the number of characters in $s$, and thus $c_i \in s$ denotes the $i^{\text{th}}$ character in $s$. Since every symbol in the alphabet of $L$ is assigned a unique natural number for the purpose of constructing the Gödel numbers and since every $c_i \in s$ represents some symbol then during the calculation of the function $f_n$ each $c_i \in s$ is calculated using the natural number that was assigned to the symbol that it represents. Finally, $\sum_{c_i \in s}$ indicates the sum ranging over every $c_i \in s$.

In definition 2, $S \in L^*$ denotes some finite ordered subset of $L$, which is some ordered sequence of 0 or more $s \in L$. Hence, $S = s_1, \ldots, s_m$, where $m \geq 0$ is the number of strings in $S$ and each string is not necessarily unique in $S$, and thus $s_i \in S$ denotes the $i^{\text{th}}$ string in $S$. Subsequently, $p_i$ denotes the $i^{\text{th}}$ prime number or $p_1 = 2$, $p_2 = 3$, $p_3 = 5$, … and the $f_n(s_i)$, which is from definition 1, is the Gödel number of $s_i \in S$. Finally, $\prod_{s_i \in S}$ indicates the product ranging over every $s_i \in S$.

Therefore, the bijective function $f_n$ is a one-to-one correspondence between every possible $s \in L$ and every natural number. However, $f_n$ is not a one-to-one correspondence between every formula in the language of some formal system and every natural number because, although every formula is some $s \in L$, not every $s \in L$ is a legitimate or well-formed formula. As a consequence, since the set of every well-formed formula in the language of some formal system is a proper subset of $L$ then the function $f_n$ merely defines an injective or one-to-one function between every well-formed formula and the natural numbers.

To make this more precise we will need to define a predicate that determines whether an arbitrary sequence of symbols in the language of some formal system is a well-formed formula. Since every formal system must have explicitly defined rules for mechanically deciding whether an arbitrary sequence of symbols in the language of the system is a well-formed formula, we can intuitively define the predicate $WFF_L$ that applies these rules for the language $L$. Thus, for any arbitrary sequence of symbols $s \in L$, $WFF_L(s)$ is true if and only if $s \in L$ is a well-formed formula. We will denote the set of every well-formed formula in the language $L$ of some formal system by $L_\varphi$ and, by employing the predicate $WFF_L$, we then define $L_\varphi$ as, $s \in L_\varphi \leftrightarrow s \in L \wedge WFF_L(s)$. As a result, if the domain of the bijective function $f_n$ is restricted to every $s \in L_\varphi$ or $f_n$: $L_\varphi \rightarrow \mathbb{N}$ then $f_n$ is an injective function, not a bijective function.

Subsequently, using the function $f_n$ from definition 1 of A Method of Constructing Gödel Numbers along with the least search operator that was defined in section 3: Notation and Terminology and the predicate $WFF_L$ that was just defined we will construct a bijective function between $L_\varphi$ and the natural numbers. First, we need to define a complete enumeration of every $s \in L_\varphi$.

**An Enumeration of Every Well-Formed Formula:** the enumeration $E_n$: $\mathbb{N} \rightarrow \mathbb{N}$ is a bijective function and $E_n(m) = [\varphi]$, where $[\varphi] = f_n(s)$ for some $s \in L_\varphi$, and thus value of the function $E_n(m)$ is the Gödel number of the $m^{\text{th}}$ well-formed formula in the language $L$ of some formal system. The function $E_n(m)$ is defined as:

$$E_n(0) = \mu k \in \mathbb{N}\ (\exists s \in L\ [f_n(s) = k \wedge WFF_L(s)])$$

$$E_n(m + 1) = \mu k \in \mathbb{N}\ (k > E_n(m) \wedge [\exists s \in L\ (f_n(s) = k \wedge WFF_L(s))])$$

We will note that if the inverse of the bijective function $f_n$, denoted by $f_n^{-1}$, is employed in the definition then the existential quantifier becomes unnecessary thus resulting in the more concise definition: $E_n(0) = \mu k \in \mathbb{N}\ (WFF_L(f_n^{-1}(k)))$ and $E_n(m + 1) = \mu k \in \mathbb{N}\ (k > E_n(m) \wedge WFF_L(f_n^{-1}(k)))$.

**A Bijection Between Every Well-Formed Formula and the Natural Numbers:** the bijective function $h_n$: $L_\varphi \rightarrow \mathbb{N}$ as well as its inverse $h_n^{-1}$: $\mathbb{N} \rightarrow L_\varphi$ are defined as:

$$h_n(s) = \exists m \in \mathbb{N}\ (E_n(m) = f_n(s))$$
$$h_n^{-1}(m) = \exists s \in L_\varphi\ (f_n(s) = E_n(m))$$

As a consequence, the definition of the bijective function $g_n$ from definition 2 of A Method of Constructing Gödel Numbers can then easily be modified to be $g_n$: $L_\varphi^* \rightarrow \mathbb{N}$, where $L_\varphi^*$ denotes the set of every finite ordered subset of $L_\varphi$, simply by replacing the $f_n(s_i)$ in the definition of $g_n$ with $h_n(s_i)$. The resulting function is a bijective function between $L_\varphi^*$ and $\mathbb{N}$. Since $L_\varphi^*$ is the set of every finite ordered sequence of well-formed formulas in the language $L$ of some formal system then $L_\varphi^*$ contains every deductive proof in the formal system.

In addition, every function that was just defined is primitive recursive because every unbounded quantifier or least search operator within the definitions of the functions can be explicitly bounded. This is because it may be easily seen that for every $s_x \in L_\varphi$ some $s_y \in L_\varphi$ can be constructed such that $WFF_L(s_y)$ and $f_n(s_x) < f_n(s_y)$ and thus also $h_n(s_x) < h_n(s_y)$. As a result, the Gödel number of this $s_y \in L_\varphi$ or $f_n(s_y)$ may then be employed as an explicit upper bound within the definitions.

## APPENDIX B: GÖDEL'S SYSTEM

We will provide the entire formal system that Gödel constructed in his proof with updated notation along with brief explanations of the components of his system. The notation was updated because much of the notation that Gödel employed in his proof is no longer in use. In all cases the updated notion is directly equivalent to the notation that Gödel employed.

First, we define the primitive symbols in the language of Gödel's formal system along with the unique natural numbers that were assigned to those symbols for the purpose of constructing the Gödel numbers of the formulas in his system. The primitive symbols are given in quotation marks.

**Primitive Symbols:** "0" = 1, "*f*" = 3, "¬" = 5, "∨" = 7, "∏" = 9, "(" = 11, ")" = 13 and then further, to the variables of type $n$, the numbers of the form $p^n$ where $p$ denotes a prime number greater than 13.

"0" represents the numeric quantity zero, "*f*" the successor operation, "¬" and "∨" are the logical operations of negation and disjunction respectively, "∏" is the universal quantifier, and then "(" and ")" are the left and right parenthesis respectively. The methodology for representing variables of a specific type $n$ is due to the simple theory of types and thus is included for technical accuracy only as mentioned in section 2.4: Gödel Numbers and the Arithmetization of Syntax.

We now define the axioms of Gödel's system. In the first group of axioms, the variables $n$ and $m$ are numeric variables and the variable $A$ is a set variable. (Technically, due to the simple theory of types, the variables $n$ and $m$ are variables of type $k$ and the variable $A$ is a variable of type $k + 1$.)

**Axioms I (Natural Number Axioms):**

**1:** $f(n) \neq 0$
**2:** $f(n) = f(m) \rightarrow n = m$
**3:** $(0 \in A \wedge \forall n\ [n \in A \rightarrow f(n) \in A]) \rightarrow \forall n\ (n \in A)$

The preceding axioms are based on the Peano Axioms and formally define the natural numbers within Gödel's system. Axiom I.1 states that zero is not the successor of any natural number. Axiom I.2 states that every natural number has a unique successor. And axiom I.3 is a set-theoretic formalization of the principle of mathematical induction. In the next group of axioms the variables $P$, $Q$ and $R$ are propositional variables.

**Axioms II (Axioms of Propositional Logic):**

**1:** $P \vee P \rightarrow P$
**2:** $P \rightarrow P \vee Q$
**3:** $P \vee Q \rightarrow Q \vee P$
**4:** $(P \rightarrow Q) \rightarrow (R \vee P \rightarrow R \vee Q)$

All of the preceding axioms are simply propositional tautologies. In the next group of axioms $x$ is a variable, $c$ is any term free for $x$ (and of the same type as $x$), $\varphi(x)$ is a formula and $\Psi$ is a formula within which the variable $x$ does not occur free.

**Axioms III (Axioms of Predicate Logic):**

**1:** $\forall x\ \varphi(x) \rightarrow \varphi(c)$

**2:** $\forall x\ (\Psi \vee \varphi(x)) \rightarrow \Psi \vee \forall x\ \varphi(x)$

Intuitively, axiom III.1 states that if the sentence $\forall x\ \varphi(x)$ is true and $c$ denotes any possible value of the variable $x$ then it may be inferred that $\varphi(c)$ is true. And Axiom III.2 simply states that if $x$ is not a free variable in the formula $\Psi$ and if, for every possible value of the variable $x$, either $\Psi$ or $\varphi(x)$ is true then it may be inferred that either $\Psi$ is true or else $\varphi(x)$ is true for every possible value of the variable $x$. In the next two axioms $A$ and $B$ are set variables and the variable $x$ is a variable whose values range over every possible member of the sets $A$ and $B$. (Technically, again due to the simple theory of types, the variable $x$ is a variable of type $k$ and the variables $A$ and $B$ are variables of type $k + 1$.)

**Axiom IV (Set-Theoretic Axiom of Comprehension):**

**1:** $\exists A\ (\forall x\ [x \in A \leftrightarrow \varphi(x)])$

This axiom is an axiom of Set Theory. The axiom states that every set is defined by some formula (or definite property) $\varphi(x)$. It is of interest to note that if the simple theory of types is not assumed then this specific formalization of the axiom has been proven to be simply inconsistent due to the discovery of the Russell paradox. However, even without the simple theory of types the modern version of this axiom, which has not been proven to be simply inconsistent, may be substituted for this axiom in Gödel's proof with essentially no change to the proof. The modern version of this axiom, called the axiom of subsets, is $\exists A\ (\forall x\ [x \in A \leftrightarrow x \in B \wedge \varphi(x)])$, where $B$ is not free in $\varphi(x)$.

**Axiom V (Set-Theoretic Axiom of Extensionality):**

**1:** $\forall x\ (x \in A \leftrightarrow x \in B) \rightarrow A = B$

This axiom is also an axiom of Set Theory. The axiom simply states that any two sets that contain exactly the same members are the same set and thus, intuitively, the axiom asserts that a set is defined by its members.

Finally, we will define the relations that Gödel employed in his proof. The purpose of the relations is to provide a formal definition of Gödel's proof-theoretic concept of provability, which is based on deductive proofs. With the exception of the updated and more descriptive notation, the relations along with their descriptions are essentially identical to how Gödel defined them in his proof.

**1:** $x|y \equiv_{df} \exists z \leq x\ (x = y \times z)$
$x$ is divisible by $y$.

**2:** $\text{IsPrime}(x) \equiv_{df} x > 1 \wedge \neg\exists z \leq x\ (z \neq 1 \wedge z \neq x \wedge x|z)$
$x$ is a prime number.

**3:** $0 \text{ PrimeOf } x \equiv_{df} 0$
$(n{+}1) \text{ PrimeOf } x \equiv_{df} \mu y < x\ (\text{IsPrime}(y) \wedge x|y \wedge y > n \text{ PrimeOf } x)$
$n$ PrimeOf $x$ is the $n^{th}$ prime number (in order of increasing magnitude) contained in $x$.

**4:** $0! \equiv_{df} 1$
$(n{+}1)! \equiv_{df} (n + 1) \times n!$
$n!$ is the factorial of $n$.

**5:** $\text{Prime}(0) \equiv_{df} 0$
$\text{Prime}(n{+}1) \equiv_{df} \mu y \leq \text{Prime}(n)! + 1\ (\text{IsPrime}(y) \wedge y > \text{Prime}(n))$
Prime($n$) is the $n^{th}$ prime number (in order of increasing magnitude).

**6:** $n \text{ TermOf } x \equiv_{df} \mu y \leq x\ (x|(n \text{ PrimeOf } x)^{y} \wedge \neg[x|(n \text{ PrimeOf } x)^{y+1}])$
$n$ TermOf $x$ is the $n^{th}$ term of the number sequence assigned to the number $x$ (for $n > 0$ and $n$ not greater than the length of this sequence).

**7:** $\text{Len}(x) \equiv_{df} \mu y \leq x\ (y \text{ PrimeOf } x > 0 \wedge y + 1 \text{ PrimeOf } x = 0)$
Len($x$) is the length of the number sequence assigned to $x$.

**8:** $x * y \equiv_{df} \mu z \leq \text{Prime}(\text{Len}(x) + \text{Len}(y))^{x+y}\ (\forall n \leq \text{Len}(x)\ [n \text{ TermOf } z = n \text{ TermOf } x] \wedge$
$\forall n \leq \text{Len}(y)\ [n > 0 \rightarrow (n + \text{Len}(x)) \text{ TermOf } z = n \text{ TermOf } y])$
$x * y$ corresponds to the operation of "concatenating" two finite number sequences.

**9:** $\text{Sym}(x) \equiv_{df} 2^{x}$
Sym($x$) corresponds to the number sequence consisting of $x$ alone (for $x > 0$).

**10:** $\mathrm{E}(x) \equiv_{df} \mathrm{Sym}(11) * x * \mathrm{Sym}(13)$
E($x$) corresponds to the operation of “enclosing within parentheses” (11 and 13 are assigned to the primitive signs “(“ and “)”, respectively).

**11:** $n \text{ Var } x \equiv_{df} n \neq 0 \wedge \exists z \leq x\ (z > 13 \wedge \mathrm{IsPrime}(z) \wedge x = z^n)$
$x$ is a VARIABLE of TYPE $n$.

**12:** $\mathrm{IsVar}(x) \equiv_{df} \exists n \leq x\ (n \text{ Var } x)$
$x$ is a VARIABLE.

**13:** $\mathrm{Neg}(x) \equiv_{df} \mathrm{Sym}(5) * \mathrm{E}(x)$
Neg($x$) is the NEGATION of $x$.

**14:** $x \text{ Dis } y \equiv_{df} \mathrm{E}(x) * \mathrm{Sym}(7) * \mathrm{E}(y)$
$x$ Dis $y$ is the DISJUNCTION of $x$ and $y$.

**15:** $x \text{ Gen } y \equiv_{df} \mathrm{Sym}(x) * \mathrm{Sym}(9) * \mathrm{E}(y)$
$x$ Gen $y$ is the GENERALIZATION of $y$ with respect to the VARIABLE $x$ (provided $x$ is a VARIABLE). Note: Gödel writes $x \prod (y)$ where we would write $\forall x\ \varphi_y(x)$.

**16:** $0 \text{ S } x \equiv_{df} x$
$(n{+}1) \text{ S } x \equiv_{df} \mathrm{Sym}(3) * (n \text{ S } x)$
$n$ S $x$ corresponds to the operation of “putting the sign ‘ $f$ ’ $n$ times in front of $x$”.

**17:** $\mathrm{Num}(n) \equiv_{df} n \text{ S } \mathrm{Sym}(1)$
Num($n$) is the NUMERAL denoting the number $n$.

**18:** $\mathrm{Type}_1(x) \equiv_{df} \exists m \leq x\ \exists n \leq x\ ([m = 1 \vee 1 \text{ Var } m] \wedge x = n \text{ S } \mathrm{Sym}(m))$
$x$ is a SIGN OF TYPE 1.

**19:** $\mathrm{Type}_n(x) \equiv_{df} (n = 1 \wedge \mathrm{Type}_1'(x)) \vee (n > 1 \wedge \exists v \leq x\ [n \text{ Var } v \wedge x = \mathrm{Sym}(v)])$
$x$ is a SIGN OF TYPE $n$.

**20:** $\mathrm{EF}(x) \equiv_{df} \exists y \leq x\ \exists z \leq x\ \exists n \leq x\ (\mathrm{Type}_n(y) \wedge \mathrm{Type}_{n+1}(z) \wedge x = z * \mathrm{E}(y))$
$x$ is an ELEMENTARY FORMULA.

**21:** $\mathrm{Op}(x, y, z) \equiv_{df} x = \mathrm{Neg}(y) \vee x = y \text{ Dis } z \vee \exists v \leq x\ (\mathrm{IsVar}(v) \wedge x = v \text{ Gen } y)$

**22:** $\mathrm{FR}(x) \equiv_{df} \mathrm{Len}(x) > 0 \wedge \forall n \leq \mathrm{Len}(x)\ (n > 0 \rightarrow [\mathrm{EF}(n \text{ TermOf } x) \vee \exists p < n\ \exists q < n\ (p > 0 \wedge$
$q > 0 \wedge \mathrm{Op}(n \text{ TermOf } x, p \text{ TermOf } x, q \text{ TermOf } x)])$
$x$ is a SEQUENCE OF FORMULAS, each term of which either is an ELEMENTARY FORMULA or results from the preceding FORMULAS through the operations of NEGATION, DISJUNCTION, or GENERALIZATION.

**23:** $\mathrm{IsFormula}(x) \equiv_{df} \exists n \leq \mathrm{Prime}([\mathrm{Len}(x)]^2)^{x \times (\mathrm{Len}(x))^2}\ (\mathrm{FR}(n) \wedge x = \mathrm{Len}(n) \text{ TermOf } n)$
$x$ is a FORMULA (that is, the last term of a FORMULA SEQUENCE $n$).

**24:** $v \text{ Bound } n, x \equiv_{df} \mathrm{IsVar}(v) \wedge \mathrm{IsFormula}(x) \wedge \exists a \leq x\ \exists b \leq x\ \exists c \leq x\ (x = a * (v \text{ Gen } b) * c \wedge$
$\mathrm{IsFormula}(b) \wedge [\mathrm{Len}(a) + 1] \leq n \leq [\mathrm{Len}(a) + \mathrm{Len}(v \text{ Gen } b)])$
the VARIABLE $v$ is BOUND in $x$ at the $n^{\text{th}}$ place.

**25:** $v \text{ Free } n, x \equiv_{df} \mathrm{IsVar}(v) \wedge \mathrm{IsFormula}(x) \wedge v = n \text{ TermOf } x \wedge n \leq \mathrm{Len}(x) \wedge \neg(v \text{ Bound } n, x)$
the VARIABLE $v$ is FREE in $x$ at the $n^{\text{th}}$ place.

**26:** $v \text{ Free } x \equiv_{df} \exists n \leq \mathrm{Len}(x)\ (v \text{ Free } n, x)$
$v$ occurs as a FREE VARIABLE in $x$.

**27:** $\mathrm{Sb}\ x\binom{n}{y} \equiv_{df} \mu z \leq \mathrm{Prime}(\mathrm{Len}(x) + \mathrm{Len}(y))^{x+y}\ (\exists u \leq x\ \exists v \leq x$
$[x = u * \mathrm{Sym}(n \text{ TermOf } x) * v \wedge z = u * y * v \wedge n = \mathrm{Len}(u) + 1])$
Sb $x\binom{n}{y}$ results from $x$ when we substitute $y$ for the $n^{\text{th}}$ term of $x$ (provided that $0 < n \leq \mathrm{Len}(x)$).

**28:** $0 \text{ St } v, x \equiv_{df} \mu n \leq \mathrm{Len}(x)\ (v \text{ Free } n, x \wedge \neg\exists p \leq \mathrm{Len}(x)\ [p > n \wedge v \text{ Free } p, x])$
$(k{+}1) \text{ St } v, x \equiv_{df} \mu n < k \text{ St } v, x\ (v \text{ Free } n, x \wedge \neg\exists p < k \text{ St } v, x\ [p > n \wedge v \text{ Free } p, x])$
$k$ St $v$, $x$ is the $(k{+}1)^{\text{th}}$ place in $x$ (counted from the right end of the FORMULA $x$) at which $v$ is FREE in $x$ (and 0 in case there is not such a place).

**29:** $\mathrm{NumFree}(v, x) \equiv_{df} \mu n \leq \mathrm{Len}(x)\ (n \text{ St } v, x = 0)$
NumFree($v$, $x$) is the number of places at which $v$ is FREE in $x$.

**30:** $\mathrm{Sub}_0(x\,{}^{v}_{y}) \equiv_{df} x$
$\mathrm{Sub}_{k+1}(x\,{}^{v}_{y}) \equiv_{df} \mathrm{Sb}\ \mathrm{Sub}_k(x\,{}^{v}_{y})\binom{k \text{ St } v, x}{y}$

**31:** $\mathrm{Sub}(x\,{}^{v}_{y}) \equiv_{df} \mathrm{Sub}_{\mathrm{NumFree}(v, x)}(x\,{}^{v}_{y})$
Sub($x\,{}^{v}_{y}$) is the notion SUBST $a\binom{v}{b}$ defined above. Note: SUBST $a\binom{v}{b}$ was defined by Gödel as the process of creating a new formula by substituting $b$ for every occurrence of the variable $v$ in the formula $a$.

**32:** a. $x \text{ Implies } y \equiv_{df} \mathrm{Neg}(x) \text{ Dis } y$
b. $x \text{ Con } y \equiv_{df} \mathrm{Neg}(\mathrm{Neg}(x) \text{ Dis } \mathrm{Neg}(y))$
c. $x \text{ Equal } y \equiv_{df} (x \text{ Implies } y) \text{ Con } (y \text{ Implies } x)$
d. $v \text{ Ex } y \equiv_{df} \mathrm{Neg}(v \text{ Gen } \mathrm{Neg}(y))$

**33:** $n \text{ Th } x \equiv_{\text{df}} \mu y \le x^{(x^n)}\ (\forall k \le \text{Len}(x)\ [(k \text{ TermOf } x \le 13 \wedge k \text{ TermOf } y = k \text{ TermOf } x) \vee (k \text{ TermOf } x > 13 \wedge k \text{ TermOf } y = k \text{ TermOf } x \times [1 \text{ PrimeOf } (k \text{ TermOf } x)^n])])$
$n$ Th $x$ is the $n^{\text{TH}}$ TYPE ELEVATION of $x$ (in case $x$ and $n$ Th $x$ are FORMULAS).
Three specific numbers, which we denote by $z_1$, $z_2$, and $z_3$, correspond to the Axioms I, 1 – 3, and we define

**34:** $\text{Z-Ax}(x) \equiv_{\text{df}} x = z_1 \vee x = z_2 \vee x = z_3$

**35:** $\text{A}_1\text{-Ax}(x) \equiv_{\text{df}} \exists y \le x\ (\text{IsFormula}(y) \wedge x = [y \text{ Dis } y] \text{ Implies } y)$
$x$ is a FORMULA resulting from Axiom Schema II, 1 by substitution. Analogously, $\text{A}_2$-Ax, $\text{A}_3$-Ax, and $\text{A}_4$-Ax are defined from Axioms [rather, Axiom Schemata] II, 2 – 4.

**36:** $\text{A-A}(x) \equiv_{\text{df}} \text{A}_1\text{-Ax}(x) \vee \text{A}_2\text{-Ax}(x) \vee \text{A}_3\text{-Ax}(x) \vee \text{A}_4\text{-Ax}(x)$

**37:** $\text{Q}(z, y, v) \equiv_{\text{df}} \neg(\exists n \le y\ \exists m \le \text{Len}(z)\ \exists w \le z\ [w = m \text{ TermOf } z \wedge w \text{ Bound } n, y \wedge v \text{ Free } n, y])$
$z$ does not contain any VARIABLE BOUND in $y$ at a place at which $v$ is FREE.

**38:** $\text{L}_1\text{-Ax}(x) \equiv_{\text{df}} \exists v \le x\ \exists y \le x\ \exists z \le x\ \exists n \le x\ (n \text{ Var } v \wedge \text{Type}_n(z) \wedge \text{IsFormula}(y) \wedge \text{Q}(z, y, v) \wedge x = [v \text{ Gen } y] \text{ Implies } [\text{Sub}(y\,{}^{v}_{z})])$
$x$ is a FORMULA resulting from Axiom Schema III, 1 by substitution.

**39:** $\text{L}_2\text{-Ax}(x) \equiv_{\text{df}} \exists v \le x\ \exists q \le x\ \exists p \le x\ (\text{IsVar}(v) \wedge \text{IsFormula}(p) \wedge v \text{ Free } p \wedge \text{IsFormula}(q) \wedge x = [v \text{ Gen } (p \text{ Dis } q)] \text{ Implies } [p \text{ Dis } (v \text{ Gen } q)])$
$x$ is a FORMULA resulting from Axiom Schema III, 2 by substitution.

**40:** $\text{R-Ax}(x) \equiv_{\text{df}} \exists u \le x\ \exists v \le x\ \exists y \le x\ \exists n \le x\ (n \text{ Var } v \wedge (n + 1) \text{ Var } u \wedge \neg[u \text{ Free } y] \wedge \text{IsFormula}(y) \wedge x = u \text{ Ex } [v \text{ Gen } ([\text{Sym}(u) * \text{E}(\text{Sym}(v))] \text{ Equal } y)])$
$x$ is a FORMULA resulting from Axiom Schema IV, 1 by substitution.
A specific number $z_4$ corresponds to Axiom V, 1, and we define

**41:** $\text{M-Ax}(x) \equiv_{\text{df}} \exists n \le x\ (n \text{ Th } z_4)$

**42:** $\text{Ax}(x) \equiv_{\text{df}} \text{Z-Ax}(x) \vee \text{A-A}(x) \vee \text{L}_1\text{-Ax}(x) \vee \text{L}_2\text{-Ax}(x) \vee \text{R-Ax}(x) \vee \text{M-Ax}(x)$
$x$ is an AXIOM.

**43:** $\text{ImmCon}(x, y, z) \equiv_{\text{df}} y = z \text{ Implies } x \vee \exists v \le x\ (\text{IsVar}(v) \wedge x = v \text{ Gen } y)$
$x$ is an IMMEDIATE CONSEQUENCE of $y$ and $z$.

**44:** $\text{ProofArray}(x) \equiv_{\text{df}} \text{Len}(x) > 0 \wedge \forall n \le \text{Len}(x)\ (n > 0 \rightarrow [\text{Ax}(n \text{ TermOf } x) \vee \exists p < n\ \exists q < n\ (p > 0 \wedge q > 0 \wedge \text{ImmCon}(n \text{ TermOf } x, p \text{ TermOf } x, q \text{ TermOf } x))])$
$x$ is a PROOF ARRAY (a finite sequence of FORMULAS, each of which is either an AXIOM or an IMMEDIATE CONSEQUENCE of two of the preceding FORMULAS).

**45:** $x \text{ ProofOf } y \equiv_{\text{df}} \text{ProofArray}(x) \wedge \text{Len}(x) \text{ TermOf } x = y$
$x$ is a PROOF of the FORMULA $y$.

**46:** $\text{Prov}(x) \equiv_{\text{df}} \exists y\ (y \text{ ProofOf } x)$
$x$ is a PROVABLE FORMULA. (Prov($x$) is the only one of the notions 1 – 46 of which we cannot assert that it is [primitive] recursive.)